\documentclass[10pt,a4paper]{amsart}
\usepackage{mathrsfs,amsthm,amssymb,amsmath,url,graphicx}

\theoremstyle{plain}
    \newtheorem{thm}{Theorem}
    
    \newtheorem{lem}[thm]{Lemma}
    \newtheorem{prop}[thm]{Proposition}
    \newtheorem{cor}[thm]{Corollary}

\theoremstyle{definition}
    \newtheorem{defn}[thm]{Definition}
    
\theoremstyle{remark}

\newcommand{\ignore}[1]{}

\newcommand{\nin}{\notin}

\DeclareMathOperator{\ran}{ran}

\newcommand{\sm}{{\setminus}}

\newcommand{\cl}[1]{\langle #1 \rangle}
 
\DeclareMathOperator{\Cl}{Cl} \DeclareMathOperator{\Cll}{Cl_{loc}}

\DeclareMathOperator{\pp}{pp}

\DeclareMathOperator{\Ideal}{Id}\DeclareMathOperator{\Mon}{Mon}

\newcommand{\nat}{\mathbb{N}}
\newcommand{\sn}{S_\nat}
\newcommand{\Q}{{\mathscr Q}}
\renewcommand{\O}{{\mathscr O}}
\newcommand{\On}{{\mathscr O}^{(n)}}

\newcommand{\Om}{{\mathscr O}^{(m)}}
\newcommand{\Oo}{{\mathscr O}^{(1)}}
\newcommand{\Ot}{{\mathscr O}^{(2)}}

\DeclareMathOperator{\Pol}{Pol} \DeclareMathOperator{\pol}{Pol}
\DeclareMathOperator{\Inv}{Inv} \DeclareMathOperator{\Aut}{Aut}

\newcommand{\inv}{^{-1}}
\newcommand{\U}{{\mathfrak U}}
\newcommand{\C}{{\mathscr C}}
\newcommand{\R}{{\mathscr R}}
\newcommand{\F}{{\mathscr F}}

\newcommand{\I}{{\mathscr I}}

\newcommand{\B}{{\mathscr B}}

\newcommand{\D}{{\mathscr D}}

\newcommand{\G}{{\mathscr G}}
\newcommand{\M}{{\mathscr M}}

\newcommand{\K}{{\mathscr K}}
\renewcommand{\H}{{\mathscr H}}

\newcommand{\Pos}{\mathbb{P}}
\renewcommand{\S}{{\mathscr S}}

\newcommand{\un}{^{(n)}}

\newcommand{\um}{^{(m)}}

\newcommand{\ut}{^{(2)}}

\newcommand{\To}{\rightarrow}
\newcommand{\mult}{\times}

\newcommand{\true}{\textit{true}}
\newcommand{\false}{\textit{false}}

\title[The Reducts of equality]
    {
        The Reducts of equality up to\\ primitive positive interdefinability
    }

\author{Manuel Bodirsky}
    \address{Laboratoire d'Informatique  (LIX), CNRS UMR 7161\\
    \`Ecole Polytechnique \\91128 Palaiseau\\
    France}
    \email{bodirsky@lix.polytechnique.fr}
    \urladdr{http://www.lix.polytechnique.fr/~bodirsky/}
\author{Hubie Chen}
    \address{Departament de Tecnologies de la Informaci\'{o} i les
    Comunicacions\\
Universitat Pompeu Fabra\\ Barcelona\\ Spain}
    \email{hubie.chen@upf.edu}
    \urladdr{http://www.tecn.upf.es/~hchen/}
\author{Michael Pinsker}
    \address{Laboratoire de Math\'{e}matiques Nicolas Oresme\\
    CNRS UMR 6139\\ Universit\'{e} de Caen\\14032 Caen Cedex\\
    France}
    \email{marula@gmx.at}
    \urladdr{http://dmg.tuwien.ac.at/pinsker/}
    \thanks{The third author is grateful for support through project P17812-N12 as well as through Erwin Schr\"{o}dinger Fellowship N2742-N18 of the Austrian Science Fund.}

\subjclass[2000]{Primary 03C40; secondary 08A40; 08A70; 03D15}

\keywords{relational structure, reduct, primitive positive
definition, lattice, invariant relation, Galois connection, local
clone, permutations}

\begin{document}

\begin{abstract}
We initiate the study of reducts of relational structures up to
primitive positive interdefinability: After providing the tools for
such a study, we apply these tools in order to obtain a
classification of the reducts of the logic of equality. It turns out
that there exists a continuum of such reducts. Equivalently,
expressed in the language of universal algebra, we classify those
locally closed clones over a countable domain which contain all
permutations of the domain.
\end{abstract}

\maketitle

{\small
    \tableofcontents
}

\section{Introduction}
Our results have impact on three fields: model theory, universal
algebra, and theoretical computer science. We therefore have a
three-fold introduction.
    \subsection{Introduction for model theorists}
In model theory, reducts of a relational structure $\Gamma$ are
usually considered \emph{up to first-order interdefinability}. To be
more precise, one considers the reducts of the expansion of $\Gamma$
by all first-order definable relations, and two reducts $\Gamma'$
and $\Gamma''$ are considered to be the same if and only if there is
a first-order definition of $\Gamma'$ in $\Gamma''$ and vice versa.
It is well-known that there is a close connection between
classifications of reducts up to first-order interdefinability and
the theory of infinite (closed) permutation groups. In 1976,
Cameron~\cite{Cameron5} showed that the highly set-transitive
permutation groups are exactly the automorphism groups of the
reducts of $(\mathbb Q,<)$. He also showed that there are exactly
five reducts of $(\mathbb Q,<)$ up to first-order interdefinability.
Recently, Junker and Ziegler gave a new proof of this fact, and
showed that $(\mathbb Q,<,a)$, the expansion of $(\mathbb Q,<)$ by a
constant $a$, has 116 reducts~\cite{JunkerZiegler}. Thomas showed
that the first-order theory of the random graph also has exactly
five reducts~\cite{RandomReducts} and conjectured that every
structure with quantifier elimination in a finite relational
signature has a finite number of reducts, up to first-order
interdefinability.

In this article, we initiate the systematic study of reducts up to
\emph{primitive positive interdefinability}. A first-order formula
is \emph{primitive positive} if and only if it is of the form
$\exists \overline x (\phi_1 \wedge \dots \wedge \phi_l)$ for atomic
formulas $\phi_1,\dots,\phi_l$.
Clearly, classifications up to primitive positive interdefinability
are harder to obtain, compared to first-order interdefinability,
since there are far more reducts to distinguish.

The simplest structure where such a classification can be studied is
the structure that has no structure at all except for the equality
relation, which is considered part of first-order logic. Our result
classifies the reducts of the structure $(X,=)$, where $X$ is a
countably infinite set, up to primitive positive definability. While
this might look trivial at first sight, we prove that the number of
such reducts equals the continuum. It turns out that the partial
order (in fact, the lattice) of these reducts is quite complex, but
can be described reasonably well. We remark that the classification
of the reducts of $(X, =)$ up to primitive positive
interdefinability is the same as the classification for a structure
of the form $(Y,=)$, where $Y$ is a set of arbitrary infinite
cardinality.

To show our classification result, we use a concept that is called
\emph{clone} in universal algebra. A \emph{clone} is a set of
operations on some fixed domain that is closed under compositions
and contains all projections. The clones we are interested in are
moreover \emph{locally closed}, which means that they are, similarly
to closed permutation groups used for the study of reducts up to
first-order interdefinability, closed sets in the natural topology
on the operations. For $\omega$-categorical structures $\Gamma$,
there is a one-to-one correspondence between the reducts of $\Gamma$
and the locally closed clones that contain the automorphisms of
$\Gamma$.
In order to classify the reducts of $(X, =)$, we thus classify those
locally closed clones which contain all permutations of the domain.

We believe that several of the results and techniques presented in
this article can be used to also classify reducts of other
$\omega$-categorical structures up to primitive positive
definability. Such future classifications are likely to combine
existing results on first-order closed reducts, such as the
above-mentioned results on the dense linear order and the random
graph, with the results of the present paper.
    \subsection{Introduction for universal algebraists}
For universal algebraists, the title of our paper could be ``The
locally closed clones above the permutations''. In fact, in this
article we describe the lattice of all locally closed clones on a
countably infinite set that contain all permutations of this set.

Traditionally, most work on clones was done for clones on a finite
domain, and only occasionally clones on infinite sets were studied.
 However recently, the use of methods from mathematical logic, in
particular from set theory, has allowed for the production of a
considerable number of new results on clones over infinite sets.
Classical results on clones, mostly on finite domains, can be found
in the books \cite{PK79}, \cite{Sze86}, and \cite{Lau06}, while the
recent survey paper \cite{GP} contains much of what is known on
clones over infinite sets.

The clones on a finite domain which contain all permutations were
completely determined in \cite{HR94}; it turns out that the number
of such clones is finite. If one moves on to infinite domains, there
are two possible clone lattices which can be considered: The lattice
of local clones, which is closer to the clone lattice on a finite
domain in both size and methodology that can be employed for its
study, and the lattice of all clones, the investigation of which
often involves set-theoretical constructions (as in \cite{GS05}, for
example).

In the full clone lattice over an infinite set, the interval of
those clones which contain all permutations has been subject to
investigation. For example, its dual atoms have been completely
described
 on domains of regular
cardinality (\cite{Hei02}, \cite{Pin05MaxAbovePerm}). Moreover, its
atoms were determined for all infinite domains in
\cite{MP06MinimalAboveS}. However, the number of clones above the
permutations has been shown to be huge (\cite{GSS-VeryMany},
\cite{Gol-AnalyticClones}, \cite{Pin05NumberOfUnary}), and it seems
unlikely that all such clones on an infinite set will ever be
classified.

As our results will show, the situation is quite different for local
clones which contain the permutations. These clones had not yet been
subject to explicit investigation and were largely unknown except
for a few implicit results, such as partial results on local clones
containing all unary operations in \cite{RS82-LocallyMaximal}. It
turns out that although the lattice of locally closed clones which
contain the permutations is uncountable, one can obtain a reasonable
understanding of it.

We remark here that the lattice of all local clones has been shown
to be quite complicated; in particular, already on a countably
infinite base set, the class of lattices which are completely
embeddable into it properly contains the class of all algebraic
lattices with a countable number of compact elements
(\cite{Pin08-local}, \cite{Pin08-morelocal}).

    \subsection{Introduction for computer scientists}
The satisfiability problem for propositional logic is one of the
most fundamental problems in computer science, and among the hardest
problems in the complexity class NP. Many restrictions of this
problem stay NP-hard -- e.g. the problems NOT-ALL-EQUAL-3SAT and
 ONE-IN-3-SAT. In a groundbreaking paper in 1978, Schaefer studied a large class of natural
restrictions of the satisfiability problem~\cite{Schaefer}. For a
fixed structure $\Gamma$ over the Boolean domain $\{0,1\}$, the
computational problem is to decide whether a given sentence of the
form $\exists \overline x. \psi_1 \wedge \dots \wedge \psi_l$ is
true, where $\psi_i$ is of the form $R(x_{i_1},\dots,x_{i_k})$ for a
$k$-ary relation $R$ in $\Gamma$. This problem is now known as
CSP$(\Gamma)$, the \emph{constraint satisfaction problem}  for the
\emph{constraint language} $\Gamma$. (It is easy to see that 3SAT,
NOT-ALL-EQUAL-3SAT, and ONE-IN-3-SAT can be modelled as
CSP$(\Gamma)$ for appropriate Boolean constraint languages
$\Gamma$.) Schaefer proved that for Boolean constraint languages
$\Gamma$ the problem CSP$(\Gamma)$ can either be solved in
polynomial time, or is NP-complete, and he gave a complete
description of the tractable and the hard Boolean constraint
languages.

Schaefer's theorem inspired many other complexity classification
projects for computational problems that are `parametrized' by a
constraint language, and similar results have been obtained not only
for constraint satisfaction, but also for other computational tasks
such as \emph{quantified constraint satisfaction}, \emph{learning of
propositional formulas}, \emph{counting solutions}, and the
\emph{maximum solution problem}. It turns out that for all these
computational problems the computational complexity does not change
if we add relations to the constraint language $\Gamma$ that have a
primitive positive definition in $\Gamma$. This implies that the
complexity of these problems only depends on the so-called
\emph{clone} consisting of all \emph{polymorphisms} of the
constraint language, i.e., the set of those operations that preserve
all relations in $\Gamma$.

The lattice of clones on a two-element set was determined already
quite some time before Schaefer, namely by Post in 1941~\cite{Post}.
And indeed, Schaefer's classification can be formulated in terms of
Post's lattice, see e.g.~\cite{Vollmer}.

Our result can be viewed as an analog of Post's lattice for the
logic of equality (rather than for propositional logic). In this
analogy, Boolean relations correspond to relations that have a
first-order definition in $(X,=)$, where $X$ is an arbitrary
infinite set. Sets of such relations have been called \emph{equality
constraint languages}~\cite{ecsps}. Hence, in this paper we study
the lattice of polymorphism clones of equality constraint languages.
Note that equality constraint languages are exactly those constraint
languages that are preserved by all permutations of the domain; in
this sense, we study those constraint languages up to primitive
positive interdefinability that have the maximal possible degree of
symmetry.

The constraint satisfaction problem for equality constraint
languages has been investigated in~\cite{ecsps}. It was shown that
CSP$(\Gamma)$ is NP-complete if the polymorphism clone of the
constraint language is the smallest element in the lattice;
otherwise, it is polynomial-time solvable. The \emph{quantified
constraint satisfaction problem (QCSP)} for equality constraint
languages $\Gamma$ has been investigated in~\cite{qecsps}.

\section{Back and forth between logic and algebra}\label{sect:backAndForth}
 We make the notions of the introduction precise, and fix some
definitions and notation from logic and universal algebra. Moreover,
we explain how to translate our classification problem from logic to
algebra and vice-versa: That is, we demonstrate in what way reducts
correspond to local clones.
   \subsection{Notions from logic}
Consider the relational structure consisting of a countably infinite
domain and the relation of equality only. For notational simplicity,
we will assume that this domain is the set $\mathbb{N}$ of natural
numbers. By a \emph{reduct} of $(\mathbb{N},=)$ we mean a structure
$\Gamma=(\mathbb{N},\R)$, where $\R$ is a set of finitary relations
on $\mathbb{N}$ all of which can be defined by a first-order formula
in the language of equality.

Let $\phi$ be a first-order formula over the relational signature
$\tau$. Then $\phi$ is called \emph{primitive positive} or
$\emph{pp}$ iff it is of the form $\exists \overline x. \psi_1
\wedge \dots \wedge \psi_l$, where $\psi_1, \dots, \psi_l$ are
atomic $\tau$-formulas. For a set $\R$ of relations on $\mathbb{N}$,
we write $\pp(\R)$ for the set of all relations on $\mathbb{N}$
which can be defined from $\R$ by pp formulas (over the signature
appropriate for $\R$).

We say that a structure $\Gamma_1=(\mathbb{N},\R_1)$ is
\emph{pp-definable} from a structure $\Gamma_2=(\nat,\R_2)$ iff
$\R_1\subseteq\pp(\R_2)$. We say that $\Gamma_1,\Gamma_2$ are
\emph{pp-interdefinable} iff they are pp-definable from each other.
The structure $\Gamma_1$ is \emph{pp-closed} iff $\R_1=\pp(\R_1)$.
The notion of pp-definability imposes a quasiorder on the reducts of
$(\nat,=)$, and the relation of pp-interdefinability is an
equivalence relation on this order. We are interested in the partial
order which is obtained by factoring the quasiorder of reducts by
this equivalence relation. That is, we consider two reducts the same
iff they are pp-interdefinable, or equivalently, we consider the
partial order of pp-closed reducts.

In fact, the pp-closed reducts form a complete lattice, where the
meet of a family of structures $(\Gamma_i)_{i\in I}$ is the
structure $(\nat,\bigcap_{i\in I}\R_i)$, and the join is
$(\nat,\pp(\bigcup_{i\in I}\R_i))$.
    \subsection{Notions from algebra}
Fix a countably infinite base set; for convenience, we take the set
$\nat$ of natural numbers. For all $n\geq 1$, denote by $\On$ the
set $\nat^{\nat^n}=\{f:\nat^n\To \nat\}$ of $n$-ary operations on
$\nat$. Then $\O=\bigcup_{n\geq 1}\On$ is the set of all finitary
operations on $\nat$. A \emph{clone} $\C$ is a subset of $\O$
satisfying the following two properties:
\begin{itemize}
    \item $\C$ contains all projections, i.e., for all $1\leq
    k\leq n$ the operation $\pi^n_k\in\On$ defined by
    $\pi^n_k(x_1,\ldots,x_n)=x_k$, and

    \item $\C$ is closed under composition, i.e., whenever
    $f\in\C$ is $n$-ary and $g_1,\ldots,g_n\in\C$ are
    $m$-ary, then the operation $f(g_1,\ldots,g_n)\in\O\um$ defined by
    $$
        (x_1,\ldots,x_m)\mapsto f(g_1(x_1,\ldots,x_m),\ldots,g_n(x_1,\ldots,x_m))
    $$
    also is an element of $\C$.
\end{itemize}

Since arbitrary intersections of clones are again clones, the set of
all clones on $\nat$, equipped with the order of inclusion, forms a
complete lattice $\Cl(\nat)$. In this paper, we are interested in
certain clones of $\Cl(\nat)$ which satisfy an additional
topological closure property: Equip $\nat$ with the discrete
topology, and $\On=\nat^{\nat^n}$ with the corresponding product
topology (Tychonoff topology), for every $n\geq 1$. A clone $\C$ is
called \emph{locally closed} or just \emph{local} iff each of its
$n$-ary fragments $\C\cap\On$ is a closed subset of $\On$.
Equivalently, a clone $\C$ is local iff it satisfies the following
interpolation property:
\begin{quote}
        For all $n\geq 1$ and all $g\in\On$, if for all finite $A\subseteq \nat^n$
        there exists an $n$-ary $f\in\C$ which agrees with
        $g$ on $A$, then $g\in\C$.
\end{quote}

Again, taking the set of all local clones on $\nat$, and ordering
them according to set-theoretical inclusion, one obtains a complete
lattice, which we denote by $\Cll(\nat)$: This is because arbitrary
intersections of clones are clones, and because arbitrary
intersections of closed sets are closed. The join of a family
$(\C_i)_{i\in I}$ is calculated as follows: First, consider the set
 of all operations on $\nat$ which can be obtained by composing
operations from $\bigcup_{i\in I}\C_i$; this set is a clone, but
might not be locally closed. For this reason, the topological
closure of this set has to be formed in addition in order to arrive
at the join in $\Cll(\nat)$. For a set of operations
$\F\subseteq\O$, we write $\cl{\F}$ for the smallest local clone
containing $\F$; if $\C$ is a local clone and $\C=\cl{\F}$, then we
say that $\C$ is \emph{generated} by $\F$. Observe that $\cl{\F}$ is
just the intersection of all local clones containing $\F$, or
equivalently the topological closure of the set of term operations
that can be built from $\F$. Using this notation, the join of the
family $(\C_i)_{i\in I}$ is simply $\cl{\bigcup_{i\in I}\C_i}$.

This paper investigates, and in some sense classifies, those local
clones on $\nat$ which contain all permutations of $\nat$. It turns
out that the number of such clones is uncountable.
    \subsection{The translation}
Let $f\in\On$ and let $\rho\subseteq \nat^m$ be a relation. We say
that $f$ \emph{preserves} $\rho$ iff $f(r_1,\ldots,r_n)\in\rho$
whenever $r_1,\ldots ,r_n\in\rho$, where $f(r_1,\ldots,r_n)$ is
calculated componentwise. This notion of preservation links finitary
relations on $\nat$ to finitary operations and is the prime tool for
studying the reducts of $(\nat,=)$, and indeed of all other
$\omega$-categorical structures, up to primitive positive
interdefinability.

For a set of relations $\R$, we write $\pol(\R)$ for the set of
those operations in $\O$ which preserve all $\rho\in\R$. The
operations in $\pol(\R)$ are called \emph{polymorphisms} of $\R$.
The following is folklore in universal algebra, see e.g.\
\cite{Sze86}.

\begin{prop}\label{prop:polRisLocalClone}
    $\pol(\R)$ is a local clone for all sets of relations
    $\R$. Moreover, every local clone is of this form.
\end{prop}

We have seen how to assign sets of operations to sets of relations;
likewise, we can go the other way. For an operation $f\in\On$ and a
relation $\rho\subseteq \nat^m$, we say that $\rho$ is
\emph{invariant} under $f$ iff $f$ preserves $\rho$. Given a set of
operations $\F\subseteq\O$, we write $\Inv(\F)$ for the set of all
relations which are invariant under all $f\in\F$. The following is
well-known.

\begin{prop}\label{prop:invFisPPClosed}
    $\Inv(\F)$ is a pp-closed set of relations for every
    $\F\subseteq\O$.
\end{prop}

Using the Galois connection defined by the operators $\pol$ and
$\Inv$ which link operations and relations, we obtain the following
well-known alternative for describing the hull operator which
assigns to a set of operations the local clone which this set
generates (confer \cite{Sze86}).

\begin{prop}\label{prop:localClosure}
    Let $\F\subseteq\O$. Then $\cl{\F}=\pol\Inv(\F)$.
\end{prop}

The reader will now expect that our notion of closure for structures
$\Gamma$, namely closure under primitive positive definitions, will
coincide with the closure operator $\Inv\Pol$ of our Galois
connection. Although this is not true for infinite structures in
general, the following theorem from \cite{BN06} states that it holds
if $\Gamma$ is not too far from the finite. (In the following, we
often identify structures $\Gamma$ with their sets of relations,
e.g., when writing $\pol(\Gamma)$.)

\begin{thm}\label{thm:ppDefinability}
    Let $\Gamma$ be $\omega$-categorical. Then
    $\Inv\pol(\Gamma)=\pp{(\Gamma)}$.
\end{thm}

Ordering locally closed permutation groups (just like clones, a
permutation group is called locally closed iff it is a closed subset
of $\O$) by set-theoretical inclusion, one obtains a complete
lattice since arbitrary intersections of local permutation groups
are local permutation groups again. Similarly, the first-order
closed structures on $\nat$ form a complete lattice. Using our
notation, we can state the well-known connection (see
e.g.~\cite{Oligo}) between locally closed permutation groups and
first-order closed reducts of $\omega$-categorical structures as
follows.

\begin{thm}\label{thm:autInv}
    Let $\Gamma$ be an $\omega$-categorical structure. The mapping
    $\G\mapsto \Inv(\G)$ sending every locally closed permutation group
    $\G$ containing the automorphism group $\Aut(\Gamma)$ of $\Gamma$ to its set of invariant
    relations is an antiisomorphism between the lattice of such
    groups and the lattice of those reducts of $\Gamma$
    which are closed under first-order definitions.
\end{thm}

Utilizing local clones, we obtain an analogous statement to Theorem
\ref{thm:autInv}  for reducts up to primitive positive
interdefinability.

\begin{thm}\label{thm:GaloisConnectionPPreductsLocalClones}
Let $\Gamma$ be an $\omega$-categorical structure. Then we have:
    \begin{enumerate}
        \item The operator $\pol$ maps every reduct of the
        first-order expansion of $\Gamma$ to a local clone
        above $\Aut(\Gamma)$.
        \item Two reducts $\Delta_1,\Delta_2$ are mapped to
        the same local clone if and only if they are equivalent with
        respect to primitive positive interdefinability.
        \item Every local clone above $\Aut(\Gamma)$ is the polymorphism clone
        of a reduct of $\Gamma$.
        \item The mappings $\pol$ and $\Inv$ are
        antiisomorphisms beween the lattice of local clones
        above $\Aut(\Gamma)$ and the lattice of those
        reducts of the first-order expansion of $\Gamma$
        which are closed under primitive positive
        definitions.
    \end{enumerate}
\end{thm}
\begin{proof}
    (1): Let $\Gamma'$ be the expansion of $\Gamma$ by all
    first-order definable relations and let $\Delta$ be a reduct of $\Gamma'$. Its is easy to see and
    well-known that $\Aut(\Gamma)=\Aut(\Gamma')$. Thus
    $\pol(\Delta)\supseteq \Aut(\Delta)\supseteq
    \Aut(\Gamma')=\Aut(\Gamma)$. By Proposition \ref{prop:polRisLocalClone}, $\pol(\Gamma)$ is a
    local clone.\\
    (2): If $\pol(\Delta_1)=\pol(\Delta_2)$, then $\pp(\Delta_1)=\pp(\Delta_2)$ by Theorem \ref{thm:ppDefinability}.
    On the other hand, if
    $\pol(\Delta_1)\neq\pol(\Delta_2)$, then
    $\cl{\pol(\Delta_1)}=\pol(\Delta_1)\neq
    \pol(\Delta_2)=\cl{\pol(\Delta_2)}$. Thus,
    $\pol\Inv\pol(\Delta_1)\neq \pol\Inv\pol(\Delta_2)$ by Proposition \ref{prop:localClosure},
    and hence $\Inv\pol(\Delta_1)\neq \Inv\pol(\Delta_2)$.
    Theorem \ref{thm:ppDefinability} then implies
    $\pp(\Delta_1)\neq\pp(\Delta_2)$.\\
    (3): Given a local clone $\C\supseteq\Aut(\Gamma)$,
    consider $\Delta:=\Inv(\C)$. By Proposition
    \ref{prop:localClosure} we have
    $\C=\cl{\C}=\pol(\Delta)$. It remains to show that
    $\Delta$ is a reduct of $\Gamma$, i.e.,  $\Delta\subseteq
    \Gamma'$, where $\Gamma'$ is the first-order expansion
    of $\Gamma$. But this is obvious, since
    $\C\supseteq\Aut(\Gamma)$ implies
    $\Delta=\Inv(\C)\subseteq\Inv\Aut(\Gamma)=\Gamma'$.\\
    (4): By (1), (2), and (3), $\pol$ is a bijective mapping from the pp-closed reducts of $\Gamma$ onto the local clones containing $\Aut(\Gamma)$.
    It is obvious from its definition that this mapping is antitone. By
    Theorem~\ref{thm:ppDefinability}, the restriction of $\Inv$ to the local
clones containing $\Aut(\Gamma)$ is the inverse of $\pol$.
\end{proof}

The above theorem tells us that classifying reducts of
$\omega$-categorical structures up to primitive positive
interdefinability really amounts to understanding the lattice of
local clones containing its automorphisms; in our case, since
clearly all permutations are automorphisms of $(\nat,=)$, we have to
investigate local clones containing all permutations of $\nat$.

\section{The result}\label{sect:result}
Via the Galois correspondence Inv--Pol, it is possible to describe
the reducts of $(\nat,=)$ either on the relational or on the
operational side. We will now give the classification on the clone
side in as much detail as is gainful at this point of the paper;
since we still want to spare the reader the technical details, we
will have to be informal at times.

An operation is called \emph{essentially unary} iff it depends on at
most one of its variables. A \emph{unary clone} is a clone all of
whose operations are essentially unary. Clearly, unary clones are
just disguised monoids of transformations on $\nat$. The first part
of our result describes the locally closed unary clones containing
the set $S_\nat$ of all permutations of $\nat$.

\begin{thm}\label{thm:unaryClones}
Let $\U$ be the lattice of locally closed unary clones containing
$S_\nat$. Then:
\begin{enumerate}
    \item $\U$ is countably infinite.
    \item $\U$ is isomorphic to the lattice of order ideals of a certain
    partial order $\mathbb{P}$ on the set of finite increasing sequences of positive natural numbers.
    \item In particular, both $\U$ and its inverse order are algebraic (that is, isomorphic to the subalgebra lattice of an algebra), and $\U$ is a distributive lattice.
    \item $\U$ is a well-quasi-order, that is, there exist
    no infinite descending chains and no infinite antichains in $\U$. (In fact, $\U$ is a better-quasi-order.)
    \item All elements of $\U$ are finitely generated over $\sn$, i.e.,
    for every $\M\in\U$ there exists a finite set $\F\subseteq\O$ such that
    $\M=\cl{\F\cup\sn}$.
\end{enumerate}
\end{thm}

Theorem~\ref{thm:unaryClones} will be proven in Section
\ref{sect:monoids}; it is also there that the partial order
$\mathbb{P}$ is defined. We remark that $\U$ is an exceptionally
well-behaved part of $\Cll(\nat)$, since $\Cll(\nat)$ is far from
satisfying any statement of the theorem:
$|\Cll(\nat)|=2^{\aleph_0}$, $\Cll(\nat)$ is not algebraic (not even
upper continuous), and $\Cll(\nat)$ does not satisfy any non-trivial
lattice identities \cite{Pin08-local}.

\begin{center}
\begin{figure}
\includegraphics[scale=0.6]{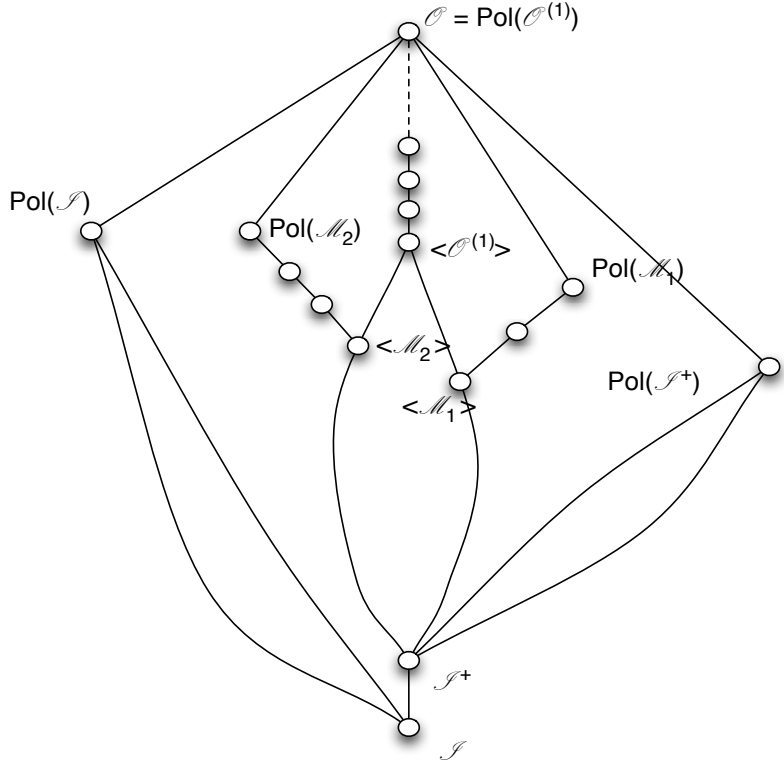}
\caption{Monoidal intervals \label{figure:MonoidalIntervals}}
\end{figure}
\end{center}

Having described the unary clones, we proceed as follows: Consider
any unary clone. Then this clone has only essentially unary
operations, and therefore differs only formally from a monoid $\M$
of transformations. Now it turns out that the set of all local
clones $\C$ which have $\M$ as their unary fragment, i.e., which
satisfy $\C\cap\Oo=\M$, forms an interval $I_\M$ of the lattice
$\Cll(\nat)$; intervals of this form are called \emph{monoidal}. The
smallest element of $I_\M$ is the unary clone we started this
argument with, namely the clone of those essentially unary
operations which are ``elements'' of $\M$; this clone is just
$\cl{\M}$. The largest element of $I_\M$ is called $\pol(\M)$ and
contains all $f\in\O$ satisfying $f(g_1,\ldots,g_n)\in\M$ whenever
$g_1,\ldots,g_n\in\M$. (This notation is consistent with our
previous use of $\pol$, if one thinks of the elements of $\M$ as
$\nat$-ary relations.) Clearly, the monoidal intervals constitute a
natural partition of $\Cll(\nat)$. Our strategy for describing the
local clones containing $S_\nat$ is to determine the monoidal
interval $I_\M$ for each locally closed monoid $\M$ containing
$S_\nat$; confer Figure \ref{figure:MonoidalIntervals}.

The following theorem describes $I_\M$ for all monoids $\M$ which
contain an operation which is neither injective nor constant. We
refer the reader to Section~5, which contains the proof of the
theorem, for the definition of quasilinearity.

\begin{thm}\label{thm:monoidalIntervalsOfLargeMonoids}
    Let $\M$ be a locally closed monoid containing $S_\nat$ as well as a non-constant and non-injective operation. Then:
    \begin{enumerate}
        \item If $\M=\Oo$, then $I_\M$ is a chain of order type
        $\omega+1$ with largest element $\O$.\\ Its smallest element is the clone $\cl{\Oo}$ of all essentially
        unary operations.\\ Its second smallest element is Burle's
        clone $\Q$ of all operations which are either essentially unary or quasilinear.\\ For $n\geq 3$, its
        $n$-th smallest element is the clone $\K_{n}$ of all
        operations which are either essentially unary or whose range
        contains less than $n$ elements.
        \item If $\M\neq \Oo$, then there exists a maximal natural
        number $k=k(\M)\geq 1$ such that $\M$ contains all unary
        operations which take at most $k$ values.\\
        If $k=1$, then $I_\M$ has only one element $\cl{\M}$.\\
        Otherwise,
         $I_\M$ is a finite chain of length $k+1$, and:\\
        Its smallest element is $\cl{\M}$.\\
        Its second smallest element $\Q^\M$ consists of $\cl{\M}$ plus all quasilinear
        operations.\\
        For $3\leq n\leq k+1$, its $n$-th smallest element $\K^\M_n$ consists of $\cl{\M}$ plus all
        operations
        whose range is smaller than $n$.
    \end{enumerate}
\end{thm}

We now turn to the monoid $\I$ locally generated by $S_\nat$. This,
as a quick check shows, consists of all injections in $\Oo$. Its
monoidal interval is the hardest to understand, and we need a few
definitions before stating the theorem describing it.

\begin{defn}[The Horn clone $\H$]\label{def:hornClone}
    Let $\H$ be the set of operations which are, up to fictitious
    variables, injective.
\end{defn}

\begin{defn}[The Bar clone $\B$]\label{defn:barClone}
    Let $f\in\Ot$ and let $k\geq 1$. If there exists an injection $p\in\Ot$ such
    that $f(x_1,x_2)=p(x_1,x_2)$ for all $x_1\geq k$ and
    $f(x_1,x_2)=p(x_1,0)=x_1$ for all $x_1<k$, then we call $f$ a
    \emph{$k$-bar function}. Let $\B$ be the clone generated (in the sense of~\ref{sect:notation}) by any
    bar function, i.e., the smallest local clone containing that bar function and all permutations of $\nat$ (we will see in Section~\ref{subsect:bar} that this definition makes sense).
\end{defn}

\begin{defn}[Richard $\R$]\label{def:injectiveInOneDirectionAndRichard}
    Let $1\leq i\leq n$. We call an operation $f\in\On$ \emph{injective in
    the $i$-th direction} iff $f(a)\neq f(b)$ whenever $a,b\in X^n$ and $a_i\neq
    b_i$. We say that $f\in\On$ is \emph{injective in one direction}
    iff there exists $1\leq i\leq n$ such that $f$ is injective
    in the $i$-th direction. Let $\R$ be the set of all operations
    which are injective in one direction.
\end{defn}

\begin{defn}[The odd clone $\S$]\label{def:oddClone}
Let $f_3 \in \O^{(3)}$ any operation satisfying the following:
    \begin{itemize}
        \item $f_3(x,1,1)= 1$, $f_3(2,x,2)= 2$, $f_3(3,3,x)= 3$, and
        \item For all other arguments, the function arbitrarily takes a value that is distinct from all other function values.
    \end{itemize}
We set $\S$ to be the clone generated by $f_3$, i.e., the smallest local clone containing $f_3$ and $S_\nat$.
\end{defn}

\begin{center}
\begin{figure}
\includegraphics[scale=0.6]{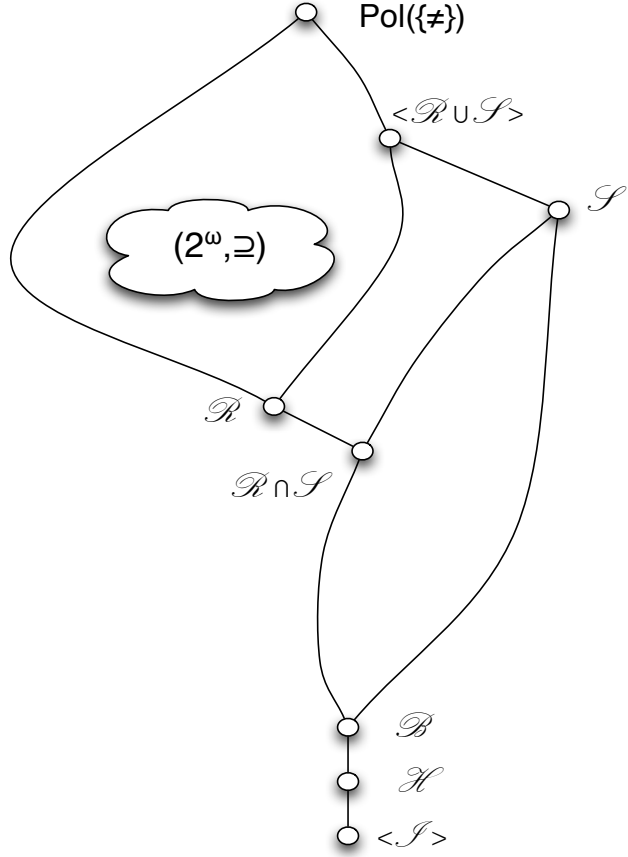}
\caption{The monoidal interval of $\I$
\label{figure:MonoidalIntervalI}}
\end{figure}
\end{center}

The following theorem summarizes the highlights of the monoidal
interval corresponding to the monoid $\I$ generated by $S_\nat$;
confer Figure~\ref{figure:MonoidalIntervalI}. More detailed
descriptions of the clones of the theorem as well as other clones in
that interval can be found in Section~\ref{sect:infiniteRangeOps}.

\begin{thm}\label{thm:monoidalIntervalOfS}
    The monoid locally generated by $S_\nat$ is the monoid $\I$ of
    injections, and:
    \begin{enumerate}
        \item The largest element of $I_\I$, $\pol(\I)$, equals
        $\pol(\{\neq\})$, where $\neq$ is the (binary) inequality
        relation.
        \item $\H$ is the unique cover of $\cl{\I}$ in $I_\I$, and all
        elements of $I_\I$ except $\cl{\I}$ contain $\H$. Moreover,
        $\H$ is generated by any binary injection, and $\Inv(\H)$
        consists of all relations definable by a Horn formula.
        \item $\B$ is the unique cover of $\H$ in $I_\I$, and all
        elements of $I_\I$ except $\cl{\I}$ and $\H$ contain $\B$.
        \item $\R,\S$ are incomparable clones in $I_\I$ and every clone in $I_\I$ is either contained in $\S$ or
        contains $\R$.
        \item The number of elements of $I_\I$ containing $\R$
        equals the continuum: In fact, the power set of $\omega$, ordered by reverse inclusion, has an order embedding into the interval $[\R,\pol(\I)]$.
        In particular, the same holds for the
        interval $I_\I$, as well as for the set of local clones above
        $S_\nat$.
    \end{enumerate}
\end{thm}

The last statement of Theorem~\ref{thm:monoidalIntervalOfS} is among
the hardest to prove in this paper, and has strong consequences for
pp classification projects, so that it deserves an own corollary.

\begin{cor}
    Let $\Gamma$ be any relational structure. Then the number of its pp-closed reducts is uncountable. In fact, there exists an
    order embedding of the power set of $\omega$ into the lattice of
    pp-closed reducts of $\Gamma$.
\end{cor}

It is for this reason that we cannot expect to completely
characterize the pp-closed reducts of any relational structure. In
our case, we obtain a complete characterization of the closed monoid
lattice and of all monoidal intervals except for those corresponding
to the monoids $\I$ and $\I^+$ (see below), where we must content
ourselves with some insights on the structure of those intervals.

It remains to describe the monoidal intervals of those monoids which
contain all injections, and some constant operations. Clearly, there
is only one such monoid, namely the monoid $\I^+$ consisting of all
constants and all injections. In general, for any set of operations
$\F\subseteq\O$, write $\F^+$ for $\F$ plus all constant operations,
and $\F^-$ for $\F$ without all constant operations. It turns out
that $I_{\I^+}$ is a complete sublattice of $I_{\I}$, as described
in the following theorem:

\begin{thm}\label{thm:monoidalIntervalOfSPlusConst}
    Let $\I^+$ be the monoid of all injective and of all constant
    operations. Then:
    \begin{enumerate}
        \item If $\C\in I_{\I^+}$, then $\C^-$ is a local clone in $I_\I$.
        \item $\cl{\I^+}^-=\cl{\I}$ and $(\pol(\I^+))^-=\S$.
        \item The mapping from $I_{\I^+}$ into the subinterval $[\cl{\I},\S]$ of $I_\I$  which sends every clone    $\C$ to $\C^-$ is a complete lattice embedding which preserves the smallest and the largest element.
        \item If $\C,\D\in I_{\I^+}$, then $\C\subsetneq \D$ iff $\C^-\subsetneq
        \D^-$.
        \item For all $\C\in I_\I$ which do not contain $\R$, $\cl{\C^+}$
        is a local clone in $I_{\I^+}$.
        \item All clones $\C$ in $I_{\I^+}$ are of this form, as  $\C=(\C^-)^+$.
        \item $\cl{\R^+}=\O$.
        \item $\H^+$ is the unique cover of $\cl{\I^+}$ in $I_{\I^+}$.
        \item $\B^+$  is the unique cover of $\H^+$  in $I_{\I^+}$.
    \end{enumerate}
\end{thm}

We remark that the mapping that sends every $\C\in I_{\I^+}$ to $\C^-$ is not surjective onto the interval $[\cl{\I},\S]$; see the remark after Proposition~\ref{prop:+notinjective}.

The following sections contain the proof of our result, and of
course more detailed definitions of the structures involved. Each
theorem corresponds to a section: Theorems~$7$, $8$, $13$ and $15$ are proven in Sections~$4$, $5$, $6$, and $7$, respectively.

    \subsection{Additional notation and terminology}\label{sect:notation}
In addition to the notation introduced so far, we establish the
following conventions. If $f\in\O$, since we are only interested in
local clones containing $S_\nat$, we abuse the notation
$\cl{\cdot}$ and write $\cl{f}$ for the local clone generated by $f$
together with $S_\nat$. For $f,g\in\O$, we say that $f$
\emph{generates} $g$ iff $g\in\cl{f}$. Similarly for $g\in \O$ and
$\F\subseteq\O$, we say that $\F$ \emph{generates} $g$ iff
$g\in\cl{\F}$. If $f\in\O$, we write $\ran(f)$ for its range.

For a relation $R$, we will usually write $\pol(R)$ instead of
$\pol(\{R\})$ for the set of all operations which preserve $R$. If
$a$ is an $n$-tuple, then we refer by $a_i$ to the $i$-th component
of $a$, for all $1\leq i\leq n$.

\section{The basis: Monoids}\label{sect:monoids}

We prove Theorem~\ref{thm:unaryClones} describing all unary clones
containing $S_\nat$. Recall that a unary clone consists only of
operations depending on at most one variable, and is therefore a
disguised monoid of transformations. For convenience, we therefore
only deal with unary operations and monoids in this section. In
particular, we adjust the meaning of certain notations for this
section: For example, $\cl{\F}$ refers to the local \emph{monoid}
(rather than the local clone) generated by a set $\F\subseteq\Oo$
together with $\sn$.

The various statements of Theorem~\ref{thm:unaryClones} are obtained
in Propositions~\ref{prop:unary:isoToIdealLattice},
\ref{prop:unary:idealLatticeWQO}, \ref{prop:unary:finitelGenerated}
and Corollaries~\ref{cor:unary:Latticedistributive} and
\ref{cor:unary:clonesWQO}.

In a first lemma, we officially state what we already observed in
the last section, namely that $\sn$ locally generates all unary
injections.

\begin{lem}
    Let $\M\supseteq\sn$ be a locally closed monoid. Then $\M$ contains the monoid $\I$ of all
    unary injective operations.
\end{lem}
\begin{proof}
    Clearly, on any finite set every injection can be locally interpolated by a
    suitable permutation.
\end{proof}

The following lemma implies that except for the full transformation
monoid $\Oo$, all closed monoids above $\sn$ consist of the
injections plus some finite range operations.

\begin{lem}\label{LEM:UNARY:noninjectiveInfRangeGeneratesAll}
    Let $f\in\Oo$ have infinite image, and assume it is not injective. Then $f$ generates all
    unary operations.
\end{lem}
\begin{proof}
    We skip the fairly easy proof, and refer the reader to the very similar (first part of the) proof of
    Lemma~\ref{LEM:UNARY:kerneltuplesorderIsGeneratingOrder}.
\end{proof}

We thus wish to know, given a finite range operation, which other
finite range operations it generates. For that, we need the
following concept.

\begin{defn}
    Let $f\in\Oo$ have finite range, and write $n=|\ran(f)|$. Enumerate the kernel classes of $f$ by $C_1,\ldots,C_n$ in such a way that their sizes
    are increasing.
    The \emph{kernel tuple} $\kappa^f\in (\omega+1)^n$ of
    $f$ is the $n$-tuple $(|C_1|,\ldots,|C_n|)$.
\end{defn}

Note that the last entry of a kernel tuple $\kappa^f$ always equals
$\omega$ since $f$ must have at least one infinite kernel class.

Having assigned a finite sequence with positive values in $\omega+1$
to every finite range operation, we are ready to order such
sequences and give the definition of $\Pos$.

\begin{defn}{\ }
    \begin{itemize}
        \item
            Let $k, n\geq 1$. For $a\in (\omega+1)^k$ and $b\in (\omega+1)^n$ we
            write
            $a\sqsubseteq b$ iff $k\leq n$ and the following holds: There
            exists a partition $\{A_1,\ldots,A_k\}$ of $\{1,\ldots,n\}$ into $k$ classes such
            that $a_i\leq \sum_{j\in A_i} b_j$ for all $1\leq i\leq k$.
        \item
            We write $\Pos_\infty$ for the partial order of finite increasing sequences of non-zero elements of $\omega+1$ ordered by $\sqsubseteq$.
        \item
            We write $\Pos$ for the partial order of the finite increasing sequences of positive natural numbers
            (not of values in $\omega+1$!) ordered by $\sqsubseteq$.
    \end{itemize}
\end{defn}

Observe that for finite increasing sequences $a,b\in (\omega+1)^n$
of the same length we have $a\sqsubseteq b$ iff $a_i\leq b_i$ for
all $1\leq i\leq n$. The following lemma justifies our definition of
$\sqsubseteq$.

\begin{lem}\label{LEM:UNARY:kerneltuplesorderIsGeneratingOrder}
    Let $f,g\in\Oo$ have finite range. Then $f$ generates $g$ iff
    $\kappa^g\sqsubseteq \kappa^f$.
\end{lem}
\begin{proof}
    Assume first that $\kappa^g\sqsubseteq \kappa^f$. Let $\kappa^f$ have length $n$ and $\kappa^g$ have length $k\leq n$. If $k=n$, then $\kappa^g_i\leq \kappa^f_i$ for all $1\leq i\leq n$. It is then not hard to see that for any finite set $A\subseteq \nat$, there exist permutations $\alpha,\beta $ such
    that $\beta\circ f\circ \alpha$ agrees with $g$ on $A$. If $k<n$, then let $A_1,\ldots,A_k$ be the partition provided by the definition
    of $\sqsubseteq$. Enumerate the kernel classes of $f$ by $C_1,\ldots,C_n$ and in such a way that $C_j$ contains $\kappa^f_j$ elements, for all
    $1\leq j\leq n$. Now take any
    $i,j\in\{1,\ldots,n\}$ which are distinct but equivalent with
    respect to the partition $A_1,\ldots,A_k$.
     By composing $f$ with a permutation, we may
    assume that $f$ maps the classes $C_i$ and $C_j$ into the class
    $C_n$, and all other classes into classes $\neq C_n$ in such a
    way that no two classes are mapped into the same class.
    Then $t_{n-1}:=f\circ f$ is a function with $n-1$ values in its range, and $\kappa^g\sqsubseteq \kappa^{t_{n-1}}$.
    Proceeding like this, we arrive after $n-k$ steps at an
    operation $t_{k}$ which takes $k$ values and which satisfies
    $\kappa^g_i\leq \kappa^{t_k}_i$ for all $1\leq i\leq k$. Thus we are back in the case $k=n$, and the proof of this direction of the lemma is finished.

    For the other direction, assume that $f$ generates $g$. Let $k,n$ be as
    before. Since the local clone generated by $f$ is the
    topological closure of the set of term operations generated by $f$, we have that
    for every $q\in\nat$, there exists a term $t_q$ consisting of permutations and $f$ which
    agrees with $g$ on the finite set $\{0,\ldots,q\}$. We can write each $t_q$ as $t_q=s_q\circ f\circ \alpha_q$, where
    $s_q$ consists of permutations and $f$, and
    $\alpha_q$ is a permutation. Thus in
    every term $t_q$, certain classes of $f$ are joined by the application of $s_q$ (and
    shifted by $\alpha_q$, which we do not care about for the moment); since there are only finitely many
    possibilities of joining classes of $f$, there is one
    constellation which appears for infinitely many $q$. Since
    $q\leq q'$ implies that $t_{q'}$ agrees with $g$ on
    $\{0,\ldots,q\}$, by replacing terms we may assume that the same classes are joined for all $q\in\nat$. Naturally,
    this partition of classes
    induces a partition $A_1,\ldots,A_w$ on $\{1,\ldots,n\}$ via
    $\kappa^f$. If $p,r\leq q$ are equivalent modulo the kernel of $g$,
    then the kernel classes of $f$ containing
    $\alpha_q(p)$ and $\alpha_q(r)$, respectively, are equivalent
    with respect to the partition $A_1,\ldots,A_w$. On the other hand, if $p,r\leq q$ are not in the same class of $g$, then $\alpha_q(p)$ and
    $\alpha_q(r)$ will not lie in equivalent $f$-classes.  Thus, taking $q$ large enough so that $\{0,\ldots,q\}$ meets all
    kernel classes of $g$, we can assign to every $g$-class (with,
    say, index $1\leq u\leq k$) an equivalence class
    $A_e\in\{A_1,\ldots,A_w\}$ in an injective way; in particular, $k\leq w$.
    For infinitely many $q$, this assignment is the same; again, by replacing terms where
    necessary, we may assume it is always the same.
    Since
    for large enough $q$ arbitrarily large parts of the kernel
    classes of $g$ are hit, we must have $\kappa^g_u\leq \sum_{j\in A_e} \kappa^f_j$ for all $1\leq u\leq k$, where $A_e$ is the class assigned to
    $u$. By joining some classes $A_i$, we can obtain $w=k$
    without changing the latter fact.
\end{proof}

\begin{lem}\label{lem:unary:unboundedFiniteRange}
 Let $\F\subseteq\Oo$. If there is no finite bound to the sizes of
 the ranges of the finite range operations in $\F$, then $\F$ generates $\Oo$.
\end{lem}
\begin{proof}
    Let $f$ be any finite range operation, and let $A\subseteq \nat $ be finite. Then there exists a finite range function $g$ which agrees with
    $f$ on $A$ and whose kernel sequence $\kappa^g$ has only one entry equal to $\omega$. Now there exists a finite range operation $h$ in $\F$ with $\kappa^g\sqsubseteq \kappa^h$, so $g$ is generated by $\F$. This proves that $f$ is generated by $\F$, and hence $\F$ generates all finite range
    operations. Clearly, any operation in $\Oo$ can be interpolated on any finite set by a finite range operation, which implies our assertion.
\end{proof}

The following is a consequence of
Lemma~\ref{LEM:UNARY:kerneltuplesorderIsGeneratingOrder}; it says
that if finitely many finite range operations join forces, the joint
generating power is not more than the sum of the generating powers
of the single operations.

\begin{prop}\label{prop:unary:Closure}
    Let $\F\subseteq\Oo$ be finite. Then $\cl{\F}=\bigcup
    \{\cl{f}:f\in\F\}$.
\end{prop}
\begin{proof}
    The non-trivial direction is to show $\cl{\F}\subseteq \bigcup
    \{\cl{f}:f\in\F\}$. If $t$ is any term made of operations in $\F$,
    then it is of the form $r\circ f$, where $r$ is a term and
    $f\in\F$. Clearly, $\kappa^t\sqsubseteq \kappa^f$, implying
    $t\in\cl{\{f\}}$. Thus, $\bigcup
    \{\cl{f}:f\in\F\}$ contains all terms that can be built from
    $\F$. This implies that the union is a monoid. Being a finite union of (topologically) closed sets $\bigcup
    \{\cl{f}:f\in\F\}$ is itself closed, and hence contains even $\cl{\F}$.
\end{proof}

We now assign ideals of $\Pos$ to local monoids containing $\sn$.

\begin{defn}
    For a local monoid $\M\supseteq\sn$, we set
    $$
    \Ideal_\infty(\M)=\{\kappa^f:  f\in\M,\, f \mbox{ has finite
        range} \}\subseteq \Pos_\infty.
    $$
    For a sequence $s\in \Pos$, write $s\ast\infty\in \Pos_\infty$ for the sequence obtained by gluing $\omega$ to the end of $s$. Now set
    $$
        \Ideal(\M)=\{s\in\Pos: s\ast\infty \in \Ideal_\infty(\M)\}.
    $$
\end{defn}

By Lemma~\ref{LEM:UNARY:kerneltuplesorderIsGeneratingOrder}, $\Ideal(\M)$ is always an ideal of $\Pos$. Conversely, we
show in the following how to get closed monoids from ideals of
$\Pos$.

\begin{defn}\label{defn:unary:limit}
    Let $k\geq 1$ and let $(s^n)_{n\in\omega}$ be an ascending sequence of $k$-tuples
    in $\Pos_\infty$. We write $\lim_n (s^n)$ for the smallest (according to $\sqsubseteq$) $k$-tuple $s$ in
    $\Pos_\infty$ satisfying $s^n\sqsubseteq s$ for all
    $n\in\omega$.\\
    For an ideal $I\subsetneq\Pos$, we let $\Mon(I)$ contain all operations
    in $\Oo$ which are either injective, or which have finite range
    and whose kernel sequence is a limit of an ascending sequence of
    $k$-tuples of the form $s\ast\infty$, where $s\in I$.
\end{defn}

\begin{lem}\label{lem:unary:monIIsMonoid}
    Let $I\subsetneq\Pos$ be an ideal. Then $\Mon(I)$ is a local monoid containing $\sn$.
\end{lem}
\begin{proof}
    Observe that if $f\in\Mon(I)$ is
    a finite range function, and if $g\in\Oo$ is a finite range function such that
    $\kappa^g\sqsubseteq \kappa^f$, then $\Mon(I)$ contains also
    $g$. For, let the range of $f$ have $n$ elements, let the range of $g$ have $k\leq n$ elements, and let $A_1,\ldots,A_k$ be the partition of $\{1,\ldots,n\}$
    provided by the definition of $\sqsubseteq$. A quick check shows that we may assume $A_k=\{n\}$. Let $(s^i)_{i\in\omega}$ be the sequence of $(n-1)$-tuples in $I$ such that $\lim_i(s^i\ast\infty)=\kappa^f$. Set $t^i_j:=\sum_{r\in A_j}s^i_r$, for all $i\in\omega$ and all $1\leq j\leq k-1$. Since $t^i\sqsubseteq s^i$, all $t^i$ are elements of $I$, provided they are actually increasing tuples. Fixing $i$, define inductively $w^i_j:=\max\{t^i_r: r\leq j\}$, for all $1\leq j\leq k-1$. It is not hard to see that the $w^i$ are in $I$, too, and that $\kappa^g$ is the limit of the increasing sequence $(w^i\ast\infty)_{i\in\omega}$, so $g\in\Mon(I)$.

    Using
    Lemma~\ref{LEM:UNARY:kerneltuplesorderIsGeneratingOrder} and
    Proposition~\ref{prop:unary:Closure},
    one now readily derives from this that $\Mon(I)$ is indeed a monoid.
    It remains to show that $\Mon(I)$ is local. Let $f\in\Oo$, and
    assume it can be interpolated on all finite sets by operations
    from $\Mon(I)$; assume also that it is not injective. If $f$ has
    infinite range, then $\Mon(I)$ must contain non-injections of
    arbitrarily large finite range, which in turn implies that $I$ contains
    tuples of arbitrary length. The definition of $\Pos$ then
    shows that $I=\Pos$, a contradiction. Thus, all non-injections
     in the local closure of $\Mon(I)$ have finite range. Assume
    again that $f$ is such a non-injection, and assume $\kappa^f$ has
    length $k$. Fix for every set $\{0,\ldots,n\}$ an operation
    $f_n\in\Mon(I)$ which agrees with $f$ on this
    set. From some $n$ on, the $f_n$ will have to take at least $k$ values, so we take the liberty of assuming that all $f_n$ have this property. An easy manipulation of the $f_n$ using the fact that $g\in\Mon(I)$ whenever $\kappa^g\sqsubseteq \kappa^{f_n}$ allows us to assume
    that every $f_n$ takes exactly $k$ values and that
    $\kappa^{f_n}\sqsubseteq\kappa^f$ for all $n\in\omega$. By
    thinning out the sequence, we may also assume that the kernel
    sequences of the $f_n$ are increasing with respect to $\sqsubseteq$. We then have
    $\lim_n(\kappa^{f_n})=\kappa^f$. Now replace each $k$-tuple
    $\kappa^{f_n}$ by a $(k-1)$-tuple $s^n\in\Pos$ with $s^n\ast\infty\sqsubseteq \kappa^{f_n}$ in such a way that $\lim_n
    (s^n\ast\infty)=\lim_n(\kappa^{f_n})$ and that the sequence $(s^n)_{n\in\omega}$ is still ascending. Clearly, $s^n\in I$ for all
    $n\in\omega$, and so $\kappa^f=\lim_n s^n$ implies $f\in\Mon(I)$.
    Hence, $\Mon(I)$ is locally closed.
\end{proof}

\begin{prop}\label{prop:unary:isoToIdealLattice}
    The mapping $\sigma: \M\mapsto \Ideal(\M)$ is an isomorphism from the
    lattice $\U$ of locally closed monoids that contain $\sn$ onto the
    lattice of ideals of $\Pos$.
\end{prop}
\begin{proof}
    That $\Ideal(\M)$ is an ideal of $\Pos$ for all monoids $\M$ follows directly from Lemma~\ref{LEM:UNARY:kerneltuplesorderIsGeneratingOrder}.
    From Lemma~\ref{lem:unary:unboundedFiniteRange} we know that
    $\Ideal(\M)=\Pos$ iff $\M=\Oo$, and obviously $\Ideal(\M)=\emptyset$
    iff $\M$ contains only injections. We have seen in Lemma~\ref{lem:unary:monIIsMonoid} that for any proper
    ideal $I$ of $\Pos$, $\Mon(I)$ is a local monoid, and a straightforward verification shows
    $\Ideal(\Mon(I))=I$, thus $\sigma$ is onto. Also, an easy check using Lemma~\ref{LEM:UNARY:kerneltuplesorderIsGeneratingOrder} shows that $\Mon(\Ideal(\M))=\M$ for every local monoid which contains $S_\nat$ and all of whose non-injections have finite range, so $\sigma$ is injective. It is obvious that both
    $\sigma$ and $\sigma\inv$
    are order-preserving.
\end{proof}

Recall that a lattice is \emph{algebraic} iff it is isomorphic to
the subalgebra lattice of an algebra.

\begin{cor}\label{cor:unary:Latticedistributive}
    The lattice $\U$ of local monoids above $\sn$ is distributive. Moreover, $\U$ and its dual order are algebraic.
\end{cor}
\begin{proof}
    By the preceding proposition, $\U$ is the lattice of ideals of a partial
    order. The assertions then follow from~\cite[p.~83]{CD73}.
\end{proof}

\begin{defn}
    A partial order is called a \emph{well-quasi-order} iff there
    are no infinite descending chains and no infinite antichains in
    it.\\
    We call a sequence $(a^n)_{n\in\omega}$ in a partial order with order relation $\leq$ \emph{bad}
    iff for no $i<j\in\omega$ we have $a^i\leq a^j$.
\end{defn}

A standard application of the infinite Ramsey's theorem shows that a
partial order is well-quasi-ordered iff it contains no bad sequence
(confer e.g.\ \cite{Die05}).

\begin{lem}\label{LEM:kernelTuplesWellFoundedNoInfAntichain}
    The set of finite sequences with values in $\omega+1$ ordered by $\sqsubseteq$ is a well-quasi-order. In particular, its suborders $\Pos_\infty$ and $\Pos$ are well-quasi-orders.
\end{lem}
\begin{proof}
    Assume that $(a^n)_{n\in\omega}$ were a bad sequence of such finite sequences. For every $n\in\omega$, let $w_n$ be the number of
    occurrences of $\omega$ in the tuple $a^n$. If the sequence $(w_n)_{n\in\omega}$ were unbounded, then we could find $n\in\omega$ such that $w_n$ is greater than
    the length of $a^0$, implying $a^0\sqsubseteq a^n$, a contradiction. Thus, we can thin out the sequence in such a way that all $w_n$ are
    equal. Let $b^n$ be the tuple obtained from $a^n$ by leaving away the $w_n$ components equal to $\omega$, for all $n\in\omega$. Clearly, the $b^n$ form
    a bad sequence of tuples with values in $\omega$.
    Now if the sequence of lengths of the $b^n$ were unbounded, then
    we could find some $i\in\omega$ such that the length of $b^{i}$ is
    greater than the sum of all components of $b^0$, hence
    $b^0\sqsubseteq b^{i}$, a contradiction. Thin out the sequence $(b^n)_{n\in\omega}$ is such a way that all tuples have the same length $k$.
    Now for all $1\leq j \leq k$,
    we thin out our sequence so that the sequence consisting of the
    $j$-th component of the $b^n$ is
    increasing; we can do this since $\omega$ is well-ordered. The remaining sequence of $b^n$ is ascending, a contradiction.
\end{proof}

In general, the ideal lattice of a given well-quasi-order need not be a
well-quasi-order. Certain well-quasi-orders which satisfy a certain strong combinatorial property and which are called \emph{better-quasi-orders}, however, do have the property that their ideal lattice is well-quasi-ordered. Although giving the definition of a better-quasi-order (see e.g.~\cite{Mil85}) would be out of scope of the present paper, we remark that it follows from the basic theory of better-quasi-orders that our well-quasi-order $\Pos$ is in fact a better-quasi-order (as Kruskal states in~\cite{Kru72}: ``All naturally known well-quasi-ordered sets which are known are better-quasi-ordered.''). Therefore, the lattice of ideals of $\Pos$ is a well-quasi-order (and, in fact, even a better-quasi-order as well). In order to spare the reader the pain of reading the definition of a better-quasi-order, we prove the following

\begin{prop}\label{prop:unary:idealLatticeWQO}
    The lattice of ideals of $\Pos$ is a well-quasi-order.
\end{prop}
\begin{proof}
    Suppose that $(I_n)_{n\in\omega}$ were a bad sequence of ideals. By taking away the first element of the sequence, we may assume that $I_n\neq \Pos$ for all $n\in\omega$. Consider for all $n\in\omega$ the set $J_n$ consisting
    of $I_n$ plus all limits of ascending chains of tuples of the same length in $I_n$ (as in Definition~\ref{defn:unary:limit}): So $J_n\subseteq\Pos_\infty$ for all $n\in\omega$. In $J_n$, every chain is bounded from above: If $C$ is a chain in $J_n$, then there is some $k\in\omega$ such that all tuples in $C$ have length at most $k$; for otherwise, $I_n$ contains sequences of arbitrary length, implying $I_n=\Pos$ contrary to our assumption. But if the elements of $C$ all have length at most $k$, then $C$ is bounded by construction of $J_n$ (i.e., adding the limits of ascending chains).

    Applying Zorn's lemma, we get that every element of
    $I_n$ is below some maximal tuple of $J_n$. By construction of $J_n$, every tuple in $\Pos$ which is below a maximal tuple of $J_n$ also is an element of $I_n$. The maximal elements
    of $J_n$ form an antichain with respect to $\sqsubseteq$.
    By Lemma~\ref{LEM:kernelTuplesWellFoundedNoInfAntichain}, the
    set of sequences in $\omega+1$, equipped with the order $\sqsubseteq$, is a well-quasi-order.
 In particular, the antichain of maximal elements of $J_n$ is finite. For every $n>0$, there exists
    a tuple in $I_0\setminus I_n$. Thus, there is a maximal tuple of $J_0$ which is not in $J_n$. We can thin out our sequence of $I_n$ so that this
    witnessing maximal tuple
    is the same for all $n>0$; denote it by $s^0$. Now we do the same for $J_1$ and all $n>1$, obtaining a maximal tuple $s^1\in J_1$ which is not in any $J_n$ with $n>1$. We continue inductively in this fashion,
    obtaining a sequence $(s^n)_{n\in\omega}$.
    By construction, this sequence is a
    bad sequence in the order of finite $(\omega+1)$-valued sequences with $\sqsubseteq$, a contradiction.
\end{proof}

\begin{cor}\label{cor:unary:clonesWQO}
    The lattice of local monoids above $\sn$ is well-quasi-ordered.
\end{cor}

\begin{prop}\label{prop:unary:finitelGenerated}
    Every local monoid $\M$ above $\sn$ is finitely generated over
    $\sn$, i.e., there exists a finite $\F\subseteq \Oo$ such that $\M=\cl{\F}$. Moreover, the number of such monoids is
    countable.
\end{prop}
\begin{proof}
    If $\M=\Oo$, then it is generated by any non-injective operation with infinite range, by Lemma~\ref{LEM:UNARY:noninjectiveInfRangeGeneratesAll}. Assume
    henceforth that $\M$ contains only injections and finite range operations.
    Set $J$ to consist of all kernel sequences of operations in $\M$. By local closure and what is by now a standard argument,
    one sees that every chain in $J$ has an
    upper bound in $J$ (confer e.g.\ the proof of
    Lemma~\ref{lem:unary:monIIsMonoid}). Hence, Zorn's lemma implies that
    every element of $J$ is below a maximal element of $J$. The maximal elements of $J$ form an antichain with respect to $\sqsubseteq$, and therefore are
    finite in number by Lemma~\ref{LEM:kernelTuplesWellFoundedNoInfAntichain}.
     Pick for
    each maximal tuple one corresponding operation in $\M$. The set $\F$ of operations thus chosen is as desired,
    by Lemma~\ref{LEM:UNARY:kerneltuplesorderIsGeneratingOrder}.

    The above argument shows that every local monoid containing $\sn$ is determined
    by a finite set of finite sequences with values in
    $\omega+1$.  There are only countably many possibilities for
    such finite sets, so the number of such monoids is countable.
\end{proof}
\section{Clones with essential finite range operations}\label{sect:finiteRangeOps}

Having understood the structure of the lattice of local monoids
containing $\sn$, we move on to describe the monoidal interval
corresponding to each such monoid. In this section, we will prove
Theorem~\ref{thm:monoidalIntervalsOfLargeMonoids}, which deals with
those monoids which contain a non-constant and non-injective
operation.

\begin{defn}
    An operation $f$ on a set $Y$ is called \emph{essential} iff it is not essentially unary, i.e., it depends on
    at least two of its variables.
\end{defn}

We will see that except for $\O$, all clones whose unary fragment
contains a non-constant and non-injective operation contain only
essential operations with finite range, and none with infinite
range; hence the title of this section.

It turns out that the monoidal intervals of the  monoids under
consideration here are all chains which can be described nicely.
This is essentially a consequence of a theorem from~\cite{HR94} for
clones on finite sets, and the power of local closure. In order to
state that theorem, we need the following definition.

\begin{defn}
    Let $Y$ be any set. We call an $n$-ary operation $f$ on $Y$ \emph{quasilinear} iff there
    exist functions $\phi_0: \{0,1\}\To Y$ and $\phi_1,\ldots,\phi_n:Y\To \{0,1\}$ such that
    $f(x_1,\ldots,x_n)=\phi_0(\phi_1(x_1)\dot{+}\ldots \dot{+}\phi_n(x_n))$, where $\dot{+}$ denotes the sum
    modulo $2$.
\end{defn}

\begin{thm}[\cite{HR94}]\label{thm:s5:Haddad-Rosenberg}
    Let $Y$ be a finite set of at least three elements and let $f$ be an essential operation on $Y$. Set $k:=|\ran(f)|$.
    Then, writing $S_Y$ for the set of all permutations on $Y$, we have:
    \begin{itemize}
        \item If $k\geq 3$, then $f$ together with $S_Y$ generate all operations on $Y$
        which take at most $k$ values.
        \item If $k=2$ and $f$ is not quasilinear, then $f$
        together with $S_Y$ generate all operations on $Y$ which take at most two
        values.
        \item If $k=2$, $f$ is quasilinear, and $|Y|$ is odd, then $f$
        together with $S_Y$ generate all quasilinear operations on $Y$.
    \end{itemize}
\end{thm}

Observe that there is no such thing as local interpolation on finite
$Y$, so ``generates'' in the theorem refers to the term closure.

The following lemma is the infinite local version of this theorem.

\begin{lem}\label{lem:s5:HaddadRosenberg:forInfinite}
    Let $f\in\On$ be essential, and assume that
    $|\ran (f)|= k$, where $2\leq k <\omega$.
    \begin{enumerate}
        \item If $k\geq 3$, then $f$
        generates all operations which take at most $k$ values.
        \item If $k=2$ and $f$ is not quasilinear, then $f$
        generates all functions which take at most $k$ values.
        \item If $k=2$ and $f$ is quasilinear, then $f$ generates
        all quasilinear operations.
    \end{enumerate}
\end{lem}
\begin{proof}
    (1): Let any operation $g$ on $|\ran g|\leq k$ be
    given. It suffices to show that for every $n\in\omega$, $f$ generates an operation $h$
    which agrees with $g$ on $\{0,\ldots,n\}$. Choosing $n$ large enough, we may
    assume that the ranges of
    both $f$ and $g$ are contained in $\{0,\ldots,n\}$. Also, again
    by making $n$ larger, we may assume that the
    restriction $f'$ of $f$ to $\{0,\ldots,n\}$ is essential
    and takes $k$ values. Then $f'$ is an essential operation on $\{0,\ldots,n\}$ which
    takes $k$ values, and hence generates all such functions by
    Theorem~\ref{thm:s5:Haddad-Rosenberg}.
    In particular, $f'$ generates a function $h'$ which agrees with
    $g$ on $\{0,\ldots,n\}$. The permutations which appear in the term which
    represents $h'$ can be extended to $X$ by the identity; occurrences of $f'$ can be replaced by $f$.
    The resulting term $h$ is a function on $X$ which still agrees
    with $g$ on $\{0,\ldots,n\}$.\\
    (2): Again, let $g$ and $n$ be given. Enlarge $n$ as before, if necessary, so that the restriction of $f$ to $\{0,\ldots,n\}$ is
    not quasilinear. Now we argue as in the preceding proof.\\
    (3): The proof works as before; $n$ only has to be chosen odd in order to
    allow application of Theorem~\ref{thm:s5:Haddad-Rosenberg}.
\end{proof}

With the preceding lemma, we see that it is quite easy to understand
what happens when we add an essential finite range operation to a
monoid. We will now show that to the monoids relevant for this
section, we in fact cannot add an essential infinite range operation
without generating $\O$. To establish this, we distinguish between
those operations which preserve the binary inequality relation
$\neq$, and those which do not. The latter case can be eliminated
right away:

\begin{prop}\label{prop:s5:essentialWithInfiniteImagePreservesNeq}
    Let $f$ be an essential operation with infinite image. Then $f$
    preserves $\neq$, or it generates all operations.
\end{prop}

The proposition will follow from the following lemma.

\begin{lem}\label{lem:s5:ConstantPaste}
    Let $f\in\On$ have infinite image, and assume it does not preserve $\neq$.
    Then $f$ generates a unary non-injective function that has infinite range.
\end{lem}

Before proving the lemma we show how the proposition follows from
it:

\begin{proof}[Proof of
Proposition~\ref{prop:s5:essentialWithInfiniteImagePreservesNeq}]
    By Lemma~\ref{lem:s5:ConstantPaste}, $f$ generates a unary
    non-injective operation with infinite range;
    Lemma~\ref{LEM:UNARY:noninjectiveInfRangeGeneratesAll} then
    implies that $f$
    generates all unary operations. Let $k \geq 1$ be arbitrary. We
    can find a unary finite range operation $g_k$ such that $g_k(f)$ is
    essential and takes at least $k$ values. By
    Lemma~\ref{lem:s5:HaddadRosenberg:forInfinite}, this implies
    that $f$ generates all operations which take not more than $k$
    values. Since $k$ was arbitrary and by local interpolation, this
    implies that $f$ generates $\O$.
\end{proof}

\begin{proof}[Proof of
Lemma~\ref{lem:s5:ConstantPaste}]
    We only have to prove something if $f$ is essential, for $f$ is itself a unary non-injective operation with infinite range otherwise.
    We also assume $f$ to depend on all of its variables.\\
    \emph{Case 1:} There exist injective functions $g_1,\ldots,g_n\in\Oo$ such that
    $f(g_1,\ldots,g_n)$ is injective. Choose $S\subseteq X$ infinite such
    that $X\sm S$ and $X\sm (g_1[S]\cup\ldots \cup g_n[S])$ are infinite. Since $f$ does not preserve $\neq$, there exist
    $c,d \in X^n$ such that $c_i\neq d_i$ for all $1\leq i\leq n$ and such that $f(c)=f(d)$. We may assume that those
    values are not in $S\cup g_1[S]\cup\ldots\cup g_n[S]$. For all $1\leq i\leq n$, set $g_i'=g_i$ on $S$ and
    $g_i'(c_1)=c_i$ and $g_i'(d_1)=d_i$. Write
    $T=S\cup\{c_1,d_1\}$. On $T$, $F=f(g_1',\ldots,g_n')$ is a function
    with infinite range which is not injective. Furthermore, the $g_i'$ can be extended to permutations
    on $X$ since they are injective and have co-infinite
    domains and ranges. This completes the first case.\\
    \emph{Case 2:} $f(g_1,\ldots,g_n)$ is not injective for all injective $g_1,\ldots, g_n \in\Oo$.
    Since $f$ is not constant,
    there exist injections $u_1,\ldots,u_n\in\Oo$ such that $f(u_1,\ldots,u_n)$ is not constant.
    By assumption for this case we have that $f(u_1,\ldots,u_n)$ is not injective.
    Thus, $f(u_1,\ldots,u_n)$
    generates a function $r\in\Oo$ which is constant except for one argument, at which it
    takes another
    value. To see this, just observe that for such an $r$, the kernel sequence $\kappa^r$ satisfies $\kappa^r\sqsubseteq\kappa^h$ for any non-constant $h\in\Oo$, and
     apply Lemma
    \ref{LEM:UNARY:kerneltuplesorderIsGeneratingOrder}.
    Consider arbitrary $g_1,\ldots,g_n\in\Oo$ such that $f(g_1,\ldots,g_n)$ is injective.
    Such functions exist since $f$ has
    infinite range. Pick an infinite $S\subseteq X$ on which
    each $g_i$ is either injective or constant. Observe that it
    is impossible that all those restrictions of $g_i$ are
    constant. Say the restriction of $g_i$ to $S$ is constant with value $c_i$ for
    $1\leq i\leq k$, where $1\leq k< n$.\\
    Since $f$ depends on its first variable, there exist
    $c_{k+1},\ldots,c_n\in X$ and $a_1,\ldots,a_k\in X$ such that
    $f(a_1,\ldots,a_k,c_{k+1},\ldots,c_n)\neq f(c)$. For otherwise
    the value of $f$ is determined by the values of the arguments
    $x_{k+1},\ldots,x_n$, contradicting that $f$ depends on $x_1$.
    Choose
    $d_{k+1},\ldots,d_n\in X$ such that
    $f(c_1,\ldots,c_k,d_{k+1},\ldots,d_n)$ is distinct from both
    $f(a_1,\ldots,a_k,c_{k+1},\ldots,c_n)$ and $f(c)$; this is
    possible as $F(x_{k+1},\ldots,x_n)=f(c_1,\ldots,c_k,x_{k+1},\ldots,x_n)$ takes
    infinitely many values. Set $u_i(0)=a_i$ and
    $u_i(x)=c_i$ if $x\neq 0$, for $1\leq i\leq k$. Let moreover
    $v_i(1)=d_i$, and $v_i(x)=c_i$ if $x\neq 1$, for all $k+1\leq i\leq n$. All $u_i$ and $v_i$ are
    generated by the two-valued function $r$ and hence by $f$.
    Set
    $F(x)=f(u_1,\ldots,u_k,v_{k+1},\ldots,v_n)(x)$. We have that
    $F(0)=f(a_1,\ldots,a_k,c_{k+1},\ldots,c_n)$ and
    $F(1)=f(c_1,\ldots,c_k,d_{k+1},\ldots,d_n)$ and $F(2)=f(c)$
    are pairwise distinct; also, $F$ has finite range. Hence,
    $f$ generates a unary function $F$ that takes exactly three
    values.\\
    Pick $a,b\in X^n$ such that $a_i\neq b_i$ for all $1\leq i\leq n$ and such that $f(a)=f(b)$.
    For all $1\leq i\leq k$, let $h_i\in\Oo$  satisfy $h_i(0)=a_i$,
    $h_i(1)=b_i$, and $h_i(x)=c_i$ for all $x\nin \{0,1\}$. Clearly,
    $h_i$ is generated by the unary function $F$ with three values, and hence also by $f$.
    For all $k+1\leq i\leq n$, let
    $\alpha_i$ be any permutation that maps $0$ to $a_i$ and $1$ to $b_i$. Now we set
    $G(x)=f(h_1(x),\ldots,h_k(x),\alpha_{k+1}(x),\ldots,\alpha_n(x))$ and have that
    $G(0)=f(a)=f(b)=G(1)$ and $G$ takes infinitely many values.
\end{proof}

We are thus left with essential infinite range operations which do
preserve the inequality relation. The crucial theorem here has been
shown in~\cite{ecsps}.

\begin{thm}\label{thm:s5:bininj}
    Every essential operation $f \in
    \pol(\neq)$ generates a binary injective operation.
\end{thm}

Able to produce a binary injection, we use this operation in order
to enlarge ranges of unary operations:

\begin{lem}\label{lem:s5:binaryInjectionPlusFiniteRange}
    Let $f\in\Ot$ be injective, and let $g\in\O$ be a non-constant
    function with finite range. Then $f$ and $g$ together generate
    $\O$.
\end{lem}
\begin{proof}
    Write $k:=|\ran(g)|$. Either $g$ is itself (essentially) unary, or it
    generates a non-constant unary operation which takes $k$ values by Lemma~\ref{lem:s5:HaddadRosenberg:forInfinite}.
    Pick such a unary operation $g_0$, and set $k_0:=k$. The
    operation $t(x,y):=f(g(x),g(y))$ takes $k_1:=(k_0)^2$
    values. Hence, $f$ and $g$ together also generate a unary operation $g_1$ which
    takes $k_1$ values. Continuing in this fashion, and since $k_0=k>1$, we obtain unary
    operations of all finite ranges, and hence also $\Oo$ by
    Lemma~\ref{LEM:UNARY:kerneltuplesorderIsGeneratingOrder}. But clearly, for every operation $h\in\Ot$, there
    exists a unary $h'\in\Oo$ such that $h(x,y)=h'(f(x,y))$. Thus,
    since $h$ was arbitrary, $f$ and $g$ generate $\Ot$, and in turn also
    $\O$, as it is well-known and easy to see that $\cl{\Ot}=\O$.
\end{proof}

We are thus ready to prove
Theorem~\ref{thm:monoidalIntervalsOfLargeMonoids}.

\begin{proof}[Proof of Theorem~\ref{thm:monoidalIntervalsOfLargeMonoids}]
    We are given a monoid $\M$ which contains an operation which is neither constant nor injective. It follows from
    Lemma~\ref{lem:s5:binaryInjectionPlusFiniteRange} that no
    clone properly contained in $\O$ and containing $\M$ can have a binary injection.
    Therefore, by Theorem~\ref{thm:s5:bininj}, it cannot have
    any essential operation which preserves the inequality relation. Nor can
    it contain an essential operation with infinite range which does
    not preserve the inequality, by
    Proposition~\ref{prop:s5:essentialWithInfiniteImagePreservesNeq}.
    Thus, it cannot contain any essential operation with infinite
    range. We now distinguish the two cases of the theorem:\\
    (1): If $\M=\Oo$, then the only clones above $\Oo$ can be the
    ones mentioned in the theorem, by
    Lemma~\ref{lem:s5:HaddadRosenberg:forInfinite}. That these
    sets of operations are actually clones is a straightforward
    verification and left to the reader.\\
    (2): If $\M\neq\Oo$, then there is a largest natural number $k$ such
    that $\M$ contains all unary operations which take at most $k$
    values; this follows from Lemma~\ref{lem:unary:unboundedFiniteRange}. If $k=1$, then no clone
    having $\M$ as its unary fragment can have an essential
    operation, since this operation would generate all quasilinear
    operations by
    Lemma~\ref{lem:s5:HaddadRosenberg:forInfinite}, and hence
    all unary operations with at most two values. Consider thus the
    case $k>1$. Again by
    Lemma~\ref{lem:s5:HaddadRosenberg:forInfinite}, the only
    clones in $I_\M$ can be $\cl{\M}$, $\Q^\M$ and the $\K_n^\M$, where $n\leq
    k+1$. It is easy to verify that these are indeed clones.
\end{proof}
\section{Clones with essential infinite range operations}\label{sect:infiniteRangeOps}
This section deals with the monoidal interval of $I_\I$, the monoid
of all injections, and contains the proof of
Theorem~\ref{thm:monoidalIntervalOfS}.
Before we start with an outline of this section, let us observe that
as a straightforward verification shows, the largest element of this
interval, $\pol(\I)$, equals $\pol(\neq)$ (this is statement (1) of
Theorem~\ref{thm:monoidalIntervalOfS}). In particular, all
operations of clones in this interval have infinite range; hence the
title of this section.

In Subsection~\ref{subsect:horn} we prove that $\H$ is the unique
cover above $\left< \I \right>$ in $I_\I$, and that $\Inv(\H)$ is
the set of all relations definable by a Horn formula
(Theorem~\ref{thm:monoidalIntervalOfS}, part (2)).

In Subsection~\ref{subsect:bar} we prove that $\B$ is the unique
cover of $\H$ in $I_\I$ (Theorem~\ref{thm:monoidalIntervalOfS}, part
(3)).

In Subsection~\ref{subsect:bar-odd} we present an infinite
strictly decreasing chain of clones that contain $\B$, and are
contained in $\S$. To this end, we also give relational descriptions
of $\B$ and of $\S$.

In Subsection~\ref{subsect:richard} we show that $\R$ and $\S$ are
incomparable, and that every clone in $I_\I$ is either contained in
$\S$ or contains $\R$ (Theorem~\ref{thm:monoidalIntervalOfS}, part
(4)).

In Subsection~\ref{subsect:richardsmanyfriends} we show that the
power set of $\omega$ embeds into the intervals of those clones in
$I_\I$ that contain $\R$; in particular, the size of this interval
equals the continuum  (Theorem~\ref{thm:monoidalIntervalOfS}, part
(5)).
    \subsection{The Horn clone}
        \label{subsect:horn}
In this subsection we show that $\cl{\I}$ has the cover $\H$ in the
monoidal interval $I_\I$.

Recall from Definition~\ref{def:hornClone} that $\H$ is the clone
consisting of all operations $f$ that are essentially injective,
i.e., all operations $f$ that are the composition $i(p_1,\dots,p_n)$
of an injective operation $i$ with projections $p_1,\dots,p_n$ (it
is straightforward to verify that this set of operations indeed
forms a locally closed clone). The clone $\H$ is also called the
\emph{Horn clone}; the reason for this name will be given in
Proposition~\ref{prop:HornChar}.

\begin{defn}
Let $\Phi$ be a quantifier-free first-order formula where all atomic
subformulas are of the form $x=y$. Then $\Phi$ is called \emph{Horn}
if $\phi$ is in conjunctive normal form (henceforth abbreviated CNF)
and each clause in $\Phi$ contains at most one positive literal.
\end{defn}

\begin{defn}
    Let $\phi(x_1,\ldots,x_n)$ be a formula in CNF. We call $\phi$ \emph{reduced} iff it is not
    logically equivalent to any of its subformulas, i.e., there is no
    formula $\psi$ obtained from $\phi$ by deleting literals or clauses
    such that $\phi(a)$ iff $\psi(a)$ for all $a \in \nat^n$.
\end{defn}

Clearly, for every formula $\phi$ in CNF there exists a logically
equivalent reduced formula.


In the following proposition, the equvalence of (1) and (7) proves
item (2) of Theorem~\ref{thm:monoidalIntervalOfS}, stating that the
Horn clone is the unique cover of $\cl{\I}$ in $I_\I$. Moreover,
items (5) and (6) provide finite relational generating systems of
$\Inv(\H)$, and (2) provides a finite operational generating system
(in fact: a continuum of such systems) of $\H$. Recall that a
formula is \emph{Horn} iff it is in conjunctive normal form and each
of its clauses contains at most one positive literal. Item (4) gives
us a syntactical description of the formulas defining relations in
$\Inv(\H)$.

\begin{prop}\label{prop:HornChar}
    For all relations $R$ with a first-order definition in $({\mathbb N}, =)$ the following are equivalent.
    \begin{enumerate}
    \item $R$ is preserved by an essential operation that preserves $\neq$.
    \item $R$ is preserved by a binary injective operation.
    \item Every reduced definition of $R$ is Horn.
    \item $R$ has a Horn definition.
    \item $R$ has a pp definition in $({\mathbb N},\neq,I)$ where
    $$I := \{(a,b,c,d) \in {\mathbb N}^4 \; | \; a=b \rightarrow c=d\}$$
    \item $R$ has a pp-definition in $({\mathbb N},N)$ where
    $$N := \{(a,b,c,d) \in {\mathbb N}^4 \; | \; a=b\neq c=d \vee \; |\{a,b,c,d\}| = 4 \}$$
    \item $R$ is preserved by $\H$.
    \end{enumerate}
\end{prop}

\begin{proof}
    The implication from (1) to (2) is Theorem~\ref{thm:s5:bininj}.

    For the implication from (2) to (3), suppose that $\Phi$ is a reduced formula that defines $R$ but is not Horn. Then there exists a clause
    $\psi$ of $\Phi$ which
    contains two equalities $x_i=x_j$ and $x_k=x_l$. Construct
    $\psi'$ from $\psi$ by removing the equation $x_i=x_j$, and
    $\psi''$ by removing $x_k=x_l$. Since $\Phi$ is reduced, there exist $a,b\in {\mathbb N}^n$ such that
    $\phi(a)$ but not $\phi'(a)$, and $\phi(b)$ but not $\phi''(b)$. Clearly, $a_i=a_j$, $a_k\neq a_l$,
    $b_i\neq b_j$, and $b_k=b_l$. Set $c=f(a,b)$, where $f$ is the binary injection preserving $R$.
    Then $c_i\neq c_j$, $c_k\neq c_l$, and in fact $\phi(c)$
    does not hold. Hence $R$ is not preserved by $f$, a
    contradiction.

    The implication from (3) to (4) is trivial.

    For the implication from (4) to (5), let $\Phi$ be a Horn formula.
    It suffices to show that all clauses $\psi$ of $\Phi$ have
    a pp definition in $({\mathbb N},I,\neq)$. If $\psi$
    is of the form $(u_1=v_1 \wedge \dots \wedge u_l = v_l) \rightarrow u=v$,
    consider the following pp formula.
    \begin{align} \exists w_1,\dots,w_{l+1}. I(w_1,w_l,u,v) \wedge
    \bigwedge_{1 \leq i \leq l} I(u_i,v_i,w_i,w_{i+1})
    \label{eq:ppdef}
    \end{align}
    Assume that $u_i=v_i$ for all $1 \leq i \leq l$. In this case, the pp formula
    implies that $w_1=w_2=\dots=w_{l+1}$, and hence also implies that $u=v$.
    Now, if $u_i \neq v_i$ for some $1 \leq i \leq l$, then for all choices
    of values for the other free variables
    the formula can be satisfied by setting $w_1,\dots,w_{l+1}$ to values   that are distinct from all other values.
    Hence, the formula is a pp definition of $\psi$ in $({\mathbb N},I)$.
    If $\psi$ does not contain a positive literal,
    consider the formula $\exists u,v. (\alpha \vee u=v) \wedge u\neq v$,
    which is equivalent to $\psi$ (we assume that $u$ and $v$ are fresh variables).
    We have seen above that the term in brackets is equivalent to a pp
    formula.
    It is then straightforward to rewrite the whole expression as a pp formula.

    For the implication from (5) to (6), it suffices to show that $I$ and $\neq$ have
    pp definitions in $({\mathbb N},N)$. For $\neq$, this is obvious.
    To express $a=b \rightarrow c=d$, consider the pp formula
    $$ \exists u,v,w. N(a,b,u,v) \wedge N(a,b,v,w) \wedge N(u,w,c,d) \; .$$
    If $a=b$, the variables $u,v,w$ must denote the same value, and hence
    the formula implies $c=d$. If $a \neq b$, then for all choices of
    values for $c$ and $d$ it is possible
    to select values for $u,v,w$ that satisfy the formula.
    Hence, the above formula is a pp definition of $I$ in $({\mathbb N},N)$.

    It is straightforward to verify the implication from (6) to (7), because
    $N$ is preserved by projections and by injective operations.

    The implication from (7) to (1) is
    immediate, since every at least binary injective operation is
    essential and has infinite image.
\end{proof}
    \subsection{The bar clone}
        \label{subsect:bar}

We show that the bar clone  $\B$ (confer
Definition~\ref{defn:barClone}) is the smallest clone below
$\pol(\neq)$ that strictly contains $\H$, thus proving item (3) of
Theorem~\ref{thm:monoidalIntervalOfS}.

A smallest non-empty $k$-ary relation that is preserved by all
permutations of $\mathbb N$ is called an \emph{orbit of $k$-tuples}.
It is clear that every relation that is preserved by all
permutations is the union of a finite number of orbits of
$k$-tuples. The following lemma will be useful here and in the
following subsections, and a proof can be found in~\cite{tcsps}.

\begin{lem}\label{lem:op-arity}
    Let $R$ be a $k$-ary relation that consists of $m$ orbits of $k$-tuples. Then every operation $f$ that violates $R$
    generates an $m$-ary operation that violates $R$.
\end{lem}

The relation $N$ has been introduced in
Proposition~\ref{prop:HornChar}; we recall it for the convenience of
the reader: $N$ is the four-ary relation $\{(a,b,c,d) \in {\mathbb
N}^4 \; | \; a=b \neq c=d \vee \; |\{a,b,c,d\}| = 4 \}$.

\begin{lem}\label{lem:binaboveHorn}
    Let $f \in \pol(\neq) \setminus \H$.
    Then $f$ generates a \emph{binary} operation that violates $N$.
    \end{lem}
\begin{proof}
    The operation $f$ must violate the relation $N$. Indeed, otherwise we would have $f\in\pol(\Inv(\H))=\H$, since all relations in $\Inv(\H)$
    have a primitive positive definition in $({\mathbb N},N)$ by Proposition~\ref{prop:HornChar} -- a
    contradiction.
    Since $N$ consists of
    two orbits of tuples, Lemma~\ref{lem:op-arity} implies that $f$ generates
    a binary operation that violates $N$.
\end{proof}

Clearly, the binary operation from Lemma~\ref{lem:binaboveHorn} is
essential and not injective.

Recall from Definition~\ref{defn:barClone} that a binary function
$f$ is called a $k$-bar function iff there exists an injection
$p\in\Ot$ such
    that $f(x_1,x_2)=p(x_1,x_2)$ for all $x_1\geq k$ and
    $f(x_1,x_2)=p(x_1,0)=x_1$ for all $x_1<k$; recall also that $\B$ is the clone generated by any bar function. We show next that $\B$ is well-defined.

Clearly, for fixed $k\geq 1$, the $k$-bar functions generate each
other. It is also easy to see that if $1\leq n\leq k$, then any
$k$-bar function generates all $n$-bar functions.

\begin{lem}\label{LEM:barFunctionsEquvalent}
    Let $k\geq 1$, and let $f\in\Ot$ be a $k$-bar function and $h\in\Ot$ be a $(k+1)$-bar function.
    Then $f$ generates $h$.
\end{lem}
\begin{proof}
    Let $f$ be a $k$-bar function, and $g$ be a variant of the
    $1$-bar which satisfies $g(k+1,x_2)=k+1$ for all $x_2\in \mathbb N$, and which is injective otherwise. Clearly, $g$ is generated
    by $f$. Set $h(x_1,x_2)=f(x_1,g(x_1,x_2))$. Let $j\leq k$.
    Then $h(j,x_2)=j$ for all $x_2\in \mathbb N$ since $f(j,z)=j$ for
    all $z\in \mathbb N$. Also, $h(k+1,x_2)=f(k+1,k+1)$ for all $x_2\in
    \mathbb N$. Moreover, if $a>k+1$ and $b,c,d\in \mathbb N$ and $(a,b)\neq (c,d)$,
    then it is easy to see that
    $h(a,b)\neq h(c,d)$. Hence, $h$ is essentially a $(k+1)$-bar
    function, the only difference to a $(k+1)$-bar function
    being the value of $h(k+1,x_2)$, which can be undone by
    composing any permutation which swaps $k+1$ and $h(k+1,0)$
    with $h$.
\end{proof}

The next lemma characterizes the binary operations in $\B$. Recall
the definition of ``injective in one direction" from Definition
\ref{def:injectiveInOneDirectionAndRichard}.

\begin{lem}\label{LEM:binaryBarFunctionsDescription}
    Let $f\in\Ot$. Then $f\in\B$ iff $f$ is injective in one direction and for all $c\in {\mathbb N}$, each of the unary operations
    $F(x)=f(x,c)$ and
    $F'(x)=f(c,x)$ is either constant or injective .
\end{lem}
\begin{proof}
    Assume that $f$ is of that form, and let $n\geq 2$ be given. We can say without
    loss of generality that $f$ is injective in
    the first direction. Write $\{0,\ldots,n\}$ as a disjoint union $A\cup B$
    in such a way that $f(c,x_2)$ is injective for all $c\in A$,
    and constant for all $c\in B$. If $|B|=j$, then there exist
    injections $\alpha, \beta$ and a $j$-bar function $g$ such that
    $\alpha (g(\beta(x_1),x_2))$ agrees with $f$ on $\{0,\ldots,n\}$. Hence, $f\in\B$.\\
    We use induction over terms to show the other direction. The
    statement is true for all bar-functions, including the
    projections. Assume it holds for $f,g_1,g_2\in\B$, and
    consider $t=f(g_1,g_2)$. Then $t$ is injective in one
    direction (this is easy to see, and also follows from Lemma~\ref{lem:RisAClone} in Subsection~\ref{subsect:richard}). Let $c\in \mathbb N$ be arbitrary and consider
    $T(x)=t(x,c)$. By induction hypothesis, we have that
    $G_1(x)=g_1(x,c)$ and $G_2(x)=g_2(x,c)$ are
    constant or injective. If both $G_1$ and $G_2$ are injective, then so is $T$,
    since $f$ is injective in one direction. If on the other hand
    both $G_1$ and $G_2$ are constant, then $T$ is constant as
    well. Assume therefore without loss of generality that $G_1$
    is injective and $G_2$ is constant with value $d$. Now if
    $F(x)=f(x,d)$ is constant, then so is $T$, and if $F$ is
    injective, then the same holds for $T$ again, and we are done.
\end{proof}

\begin{lem}\label{LEM:producingInjectiveFields}
    Let $f\in\pol(\neq)\ut$ be essential
    and non-injective, and let $k\geq 2$.  Then $f$ generates a
    binary function $g$ which is not injective but injective on $k=\{0,\ldots,k-1\}$.
\end{lem}
\begin{proof}
    Without loss of generality we assume $f(c,d)=f(c',d)$ for some
    $d\in \mathbb N$ and distinct $c, c'\in \mathbb N$. Pick any binary injection $h$.
    By Lemma \ref{prop:HornChar}, $h$ is generated by $f$.
    Set $f'(x_1,x_2)=h(f(x_1,x_2),x_2)$. Then $f'$
    satisfies $f'(c,d)=f'(c',d)$, and is injective in the $2$-nd
    direction. Therefore, replacing $f$ by $f'$, we may assume that $f$
    is injective in the $2$-nd direction. We use induction over $k$ to prove the assertion of the lemma.\\
    For the induction beginning, let $k=2$. We may assume that there
    exist $s>k$ and distinct $t, t'>k$ with
    $f(t,s)=f(t',s)$. Indeed, otherwise $f$ is injective on
    the infinite square of pairs from $\{k+1,\ldots\}$, in which case we get the assertion by applying permutations. Since $f$ is essential, we may assume $f(0,0)\neq f(1,0)$. We have that
    $f(0,1)\nin\{f(0,0),f(1,0)\}$ and $f(1,1)\nin\{f(0,0),f(1,0)\}$ as
    $f$ is injective in the $2$-nd direction.
    Therefore, the only thing
    in our way to an injection on $\{0,1\}$ is that possibly $f(0,1)=f(1,1)$.
    In that case, we let a permutation $\alpha$ map $f(0,1)$ to $0$, $f(t,s)$ to $s$, and set
    $f'(x_1,x_2)=f(x_1,\alpha (f(x_1,x_2)))$.
    It is easy to check that $f'$ is injective on $2$ and
    satisfies $f(t,s)=f(t',s)$.\\
    For the induction step, let $f$ be injective on $k$, where $k\geq 2$.
    We generate a function $g$ that is injective on
    $k+1$ and satisfies $g(t,s)=g(t',s)$ for some $s,t,t'\in \mathbb N$ with $t\neq t'$.
    Again, we pick a binary injection $h$ and set
    $f'(x_1,x_2)=h(f(x_1,x_2),x_2)$. Clearly, $f'$ is still injective on $k$. Moreover, $f'$
    satisfies $f'(c,d)=f'(c',d)$, and is injective in the $2$-nd
    direction. Therefore, replacing $f$ by $f'$, we may henceforth assume that $f$
    is injective in the $2$-nd direction.\\
    As in the induction beginning, we may assume that
    there exist $s>k$ and distinct $t, t'>k$ with $f(t,s)=f(t',s)$. Since $f$ is injective in
    the $2$-nd direction we have $f(t,s)\nin f[(k+1)^2]$.\\
    If $f(k,n)=f(i,n)$ for some $i< k$ and $n\leq k$, then we do the following: Choose a
    permutation $\alpha$ that exchanges $k$ with some $j< k$ for which $j\neq i$,
    and which
    is the identity otherwise. Let $\beta$ be a permutation that maps
    $\{f(0,n),\ldots,f(k,n)\}$ into $k$, and $f(t,s)$ to $s$. Now set
    $f'(x_1,x_2)=f(\alpha(x_1),\beta (f(x_1,x_2)))$. We have:
    \begin{itemize}
        \item All points $u,v\in {\mathbb N}^2$ that had different values under $f$ still
        have different values under $f'$; in particular, $f'$ is injective in
        the $2$-nd direction, and injective on $k$.
        \item $f'(i,n)\neq f'(k,n)$.
        \item $f'(t,s)=f'(t',s)$.
    \end{itemize}
    If we repeat this procedure for all $n\leq k$ for which
    there exists $i< k$ with $f(k,n)=f(i,n)$, then in the end $f(j,n)\neq f(k,n)$ for all $j<k$ and all $n\leq k$.
    Now the only non-injectivity that could be left on $k+1$ could be that
    $f(i,k)=f(j,k)$ for some distinct $i,j< k$. One gets rid of this setting a last time
    $g(x_1,x_2)=f(x_1,\beta (f(x_1,x_2)))$, where $\beta$ is a permutation that maps
    $\{f(0,k),\ldots,f(k,k)\}$ into $k$, and $f(t,s)$ to
    $s$.
\end{proof}

\begin{lem}\label{LEM:EveryBinaryGeneratesTheOneBar}
    Let $f\in\pol(\neq)\ut$ be essential and non-injective. Then it generates a $1$-bar function.
\end{lem}
\begin{proof}
    Fix a $1$-bar function $b$.
    We show that for arbitrary $k\geq 2$, $f$ generates a function which agrees with
    $b$ on $k$. Define a finite sequence of natural numbers $(m_j)_{j\leq k}$
    by setting  $m_0=k$, and $m_{j+1}=(m_j)^2$. By
    Lemma \ref{LEM:producingInjectiveFields} and the use of permutations, we can produce a function
    $g$ which is injective on $\{1,\ldots,m_{k}-1\}\mult m_{k}$ and satisfies $g(0,c)=g(0,c')=c$ for
    some distinct $c,c'\in \mathbb N$. Moreover, we may assume that $g[\{1,\ldots,m_{k-1}-1\}\mult m_{k-1}]\subseteq m_k$.
    We first want to produce a function $g'$
    which is injective on $\{1,\ldots,m_{k-1}-1\}\mult m_{k-1}$ and satisfies
    $g'(0,0)=g'(0,c)=g'(0,c')=c$. There is nothing to show if
    $g(0,0)=c$. Otherwise, we may assume $g(0,0)=c'$, and set
    $g'(x_1,x_2)=g(x_1,g(x_1,x_2))$; it is easy to check that
    $g'$ has the desired properties. Replace $g$ by $g'$. By repeating this
    procedure we can produce a function $g'$ which is injective on $\{1,\ldots,m_{k-2}-1\}\mult m_{k-2}$
    and satisfies
    $g'(0,1)=g'(0,0)=g'(0,c)=g'(0,c')=c$.
    Doing the same $k$
    times, we get $g'$ which is injective on $\{1,\ldots,m_0-1\}\mult
    m_0$ and satisfies $g'(0,k-1)=\ldots=g'(0,0)$. Clearly, there
    exists a permutation $\alpha$ such that $\alpha (g')$
    agrees with $b$ on $k^2$.
\end{proof}

The following summarizes the results about $\B$ obtained so far.
Relational descriptions of $\B$ will be given in
Proposition~\ref{prop:RelBarChar}. Observe that the equivalence
between (1) and (5) is exactly statement (3) of
Theorem~\ref{thm:monoidalIntervalOfS}.

\begin{prop}\label{prop:BarChar}
    Let $R$ be a relation with a first-order definition in $({\mathbb N}, =)$. Then the following are equivalent.
    \begin{enumerate}
    \item $R$ is preserved by an operation from $\pol(\neq)\sm\H$, i.e., $\pol(R) \cap \pol(\neq) \not\subseteq \H$.
    \item $R$ is preserved by a non-injective operation
    that depends on all its arguments and preserves $\neq$.
    \item $R$ is preserved by an operation from $\pol(\neq)^{(2)}$
      that violates $N$.
    \item $R$ is preserved by a $1$-bar function.
    \item $R$ is preserved by $\B$, i.e., $\B \subseteq \pol(R)$.
    \end{enumerate}
\end{prop}

\begin{proof}
    Given an operation from $\pol(\neq)\sm\H$, we obtain an
    operation as described in (2) by leaving away all fictitious
    arguments, so (1) implies (2).
    An operation that is non-injective and depends on all its arguments
    is not generated by an injective operation, and therefore
    Lemma~\ref{lem:binaboveHorn}
    shows the implication from $(1)$ to $(2)$. Any binary operation that
    violates $N$ is necessarily essential and non-injective, and hence
    the conditions of Lemma~\ref{LEM:EveryBinaryGeneratesTheOneBar} are
    satisfied, which shows the implication from $(2)$ to $(3)$.
    The equivalence of $(3)$ and $(4)$ is Lemma~\ref{LEM:barFunctionsEquvalent},
    and the implication from $(3)$ to $(1)$ is immediate.
    The equivalence of $(4)$ and $(5)$ is by definition.
\end{proof}
    \subsection{From the bar clone to the odd clone}
        \label{subsect:bar-odd}
In this subsection, we explore those clones in the monoidal interval
$I_\I$ which contain the bar clone $\B$, but which are not contained
in $\R$ (the clone of operations which are injective in one
direction,
cf.~Definition~\ref{def:injectiveInOneDirectionAndRichard}). It will
be necessary to first give a relational description of $\B$.

\subsubsection{The bar clone, relationally}
In the following it will be convenient to work with first-order
formulas whose atomic formulas are of the form $x=y$, \true, or
\false. The \emph{graph} of a quantifier-free formula $\Phi$ with
such atomic formulas is the graph where the vertices are the
variables of $\Phi$, and where two vertices $x,y$ are adjacent iff
$\Phi$ contains the sub-formula $x=y$ or the sub-formula $x \neq y$.
We recall standard terminology. If a formula $\Phi$ is in
conjunctive normal form, the conjuncts in $\Phi$ are called
\emph{clauses}, and the disjuncts in a clause are called
\emph{literals}. Hence, literals are formulas that are either
atomic, in which case they are also called \emph{positive}, or the
negation of an atomic formula, in which case they are also called
\emph{negative}.

\begin{defn}
Let $\Phi$ be a quantifier-free first-order formula where all atomic
formulas are of the form $x=y$ or \false. Then $\Phi$ is
called
\begin{itemize}
\item \emph{Horn} iff $\Phi$ is in conjunctive normal form and
each clause in $\Phi$ contains at most one positive literal;
\item \emph{connected Horn} iff $\Phi$ is Horn and
the graph of a clause $\psi$ in $\Phi$ has at most two connected
components when $\psi$ has no positive literals, and is connected
when $\psi$ has positive literals.
\item \emph{extended Horn} iff $\Phi$ is a conjunction of formulas of the form
$(\phi_1 \wedge \dots \wedge \phi_l) \rightarrow (\psi_1 \wedge
\dots \wedge \psi_k)$ for $l \geq 1, k \geq 1$
where $\phi_1,\dots,\phi_l,\psi_1, \dots, \psi_k$ are atomic
formulas (and hence of the form $x=y$ or \false).
\item \emph{connected extended Horn} iff $\Phi$ is extended Horn,
and if a conjunct $\psi$ of $\Phi$ is connected whenever the
right-hand side of the implication in $\psi$ contains a literal of
the form $x=y$ and no literal \false.
\end{itemize}
\end{defn}

Clearly, every connected Horn formula can be written as a connected extended
Horn formula. But there are connected extended Horn formulas that
are not equivalent to any connected Horn formula; we will even see
that there are connected extended Horn formulas that do not have a
pp definition by connected Horn formulas
(Lemmas~\ref{lem:gkpresrk-1} and \ref{lem:gkviolatesrk}).

The conjuncts in an extended Horn formula are also called
\emph{extended Horn clauses}, and the literals on the left-hand side
(right-hand side) of the implication of an extended Horn clause are
called the \emph{negative literals} (\emph{positive literals},
respectively) of the extended Horn clause. Note that extended Horn
formulas can always be translated into (standard) Horn formulas: if
$(\phi_1 \wedge \dots \wedge \phi_k) \rightarrow (\psi_1 \wedge \dots \wedge \psi_k)$ is an
extended Horn clause, then we can replace this clause by the
conjunction of Horn clauses $(\neg \phi_1 \vee \dots \vee \neg \phi_l \vee \psi_1) \wedge \dots
\wedge (\neg \phi_1 \vee \dots \vee \neg \phi_l \vee \psi_k)$.


It will be convenient to say that a formula is \emph{preserved} by
an operation (or a set of operations) iff the relation that is
defined by the formula is preserved by the operation (or by the set
of operations, respectively). We want to give a syntactic
description of the formulas that are preserved by $\B$. We show that
every such formula is equivalent to a connected extended Horn
formula. In fact, we show the stronger result that every
\emph{expanded Horn formula} preserved by $\B$~is itself a connected
extended Horn formula.

\begin{defn}\label{def:expand}
A extended Horn formula $\Phi$ is called \emph{expanded Horn} iff it
contains all connected extended Horn clauses on the same set of
variables as $\Phi$ that are implied by $\Phi$, and for every
disconnected clause $\psi_1 \rightarrow \psi_2$ in $\Phi$ where
$\psi_2$ does not contain \false\ and for all variables $x,y$ in
$\Phi$
\begin{itemize}
\item adding an atomic formula of the form $x=y$ to $\psi_1$, or
\item adding an atomic formula of the form $x=y$ to $\psi_2$, or
\item setting $\psi_2$ to \false
\end{itemize}
results in a formula that is not equivalent to $\Phi$.
\end{defn}

\begin{lem}\label{lem:expand}
Every extended Horn formula is equivalent to an expanded Horn
formula.
\end{lem}
\begin{proof}
Let $\Phi$ be any given extended Horn formula. We construct the
expanded Horn formula  that is equivalent to $\Phi$, as follows. We
first add to $\Phi$ all connected extended Horn clauses that are
implied by $\Phi$.

Then, if $\Phi$ contains an extended Horn clause $\psi_1 \rightarrow
\psi_2$ that is not connected and where $\psi_2$ does not contain
\false, and $\Phi$ is equivalent to a formula $\Phi'$ obtained from
$\Phi$ by adding an atomic formula of the form $x=y$ to $\psi_1$ or
to $\psi_2$ for variables $x,y$ in $\Phi$, or by setting $\psi_2$ to
\false, then we replace $\Phi$ by $\Phi'$. If $\Phi$ does not
contain such a clause (which will always happen after a finite
number of steps, because there is only a finite number of formulas
of the form $x=y$ for variables $x,y$ in $\Phi$), then $\Phi$ is
expanded Horn, and clearly equivalent to the formula we started
with.
\end{proof}

\begin{lem}\label{lem:connextHorn}
If $\Phi$ is an expanded Horn formula, and $\Phi$ defines a relation
$R$ that is preserved by $\B$, then $\Phi$ is connected extended
Horn.
\end{lem}
\begin{proof}
Suppose for contradiction that $\Phi(x_1,\dots,x_n)$ contains an
extended Horn clause $\psi$ of the form $\psi_1 \rightarrow \psi_2$
whose graph $G$ contains at least two connected components, and
where $\psi_2$ does not contain \false\ and contains at least one
positive literal $x_i=x_j$. Let $x_k$ be a variable from a connected
component of $G$ that does not contain $x_i$ and $x_j$. Let
$C_1,\dots,C_p$ be the components of $G' := G - H$, where $H$ is the
graph of $\psi_2$. Observe that $x_i$ and $x_j$ are in distinct
components: Otherwise, if we replace the clause $\psi$ by $(\psi_1
\wedge x_i=x_k) \rightarrow \psi_2$, we obtain an equivalent
formula. Because $\psi$ is disconnected, this is in contradiction to
the assumption that $\Phi$ is expanded. Assume for the sake of
notation that $C_1$ is the component of $x_i$, and $C_2$ the one of
$x_j$.

We claim that there is a tuple $s$ for which $\Phi(s)$ is true, for
which $s_i \neq s_j$, and for which $s_l = s_m$ whenever $x_l, x_m$
are both in $C_1$ or are both in $C_2$. Otherwise, $\Phi$ implies
$(\bigwedge_{(y,z) \in C_1^2 \cup C_2^2} y=z)
 \rightarrow x_i = x_j$, which is a connected expanded Horn formula
 and therefore contained in $\Phi$.
Then $\Phi$ must be equivalent to the formula $\Phi'$ obtained from
$\Phi$ by replacing the disconnected clause $\psi$ by $(\psi_1
\wedge x_i=x_j) \rightarrow \psi_2$, again in contradiction to the
assumption that $\Phi$ is expanded. To show the equivalence, it
clearly suffices to prove that $\Phi'$ implies $\Phi$. So suppose
that $s$ satisfies $\Phi'$ and $\psi_1$. The clause
$(\bigwedge_{(y,z) \in C_1^2 \cup C_2^2} y=z)
 \rightarrow x_i = x_j$, which is contained in $\Phi$ and in $\Phi'$,
 shows that $s_i=s_j$. But then the premise of the new extended
 Horn clause of $\Phi'$ is satisfied as well, and therefore $s$
 satisfies $\psi_2$, which is what we had to show.

Next, we claim that there is a tuple $t$ with $t_k \neq t_i$
satisfying $\Phi$ and where $t_l = t_m$ if $x_l, x_m$ are from the
same component $C_q$ for some $1\leq q \leq p$. Otherwise, $\Phi$
implies $\psi_1 \rightarrow x_k = x_i$. But then the formula $\Phi'$
obtained from $\Phi$ by replacing $\psi$ by $\psi_1 \rightarrow
(\psi_2 \wedge x_k = x_i)$ is clearly equivalent to $\Phi$. This
again contradicts the assumption that $\Phi$ is expanded.

Let $f$ be a binary operation such that $x \mapsto f(x,v)$ is
constant for all entries $v$ of $t$ except for $v=t_i$, and which is
injective otherwise. Clearly, modulo permutations acting on its
arguments, $f$ is an $n$-bar operation, and hence  $f \in \B$. Now,
consider the tuple $r := f(s,t)$. Note that $t_i=t_j$, because $t$
satisfies $\Phi$ and $t$ satisfies the premise of the clause $\psi$
that is contained in $\Phi$; since the conclusion of $\psi$ contains
$x_i=x_j$ we have that $t_i=t_j$. It is straightforward to verify
that $r$ does not satisfy $\Phi$, because $r$ satisfies $\psi_1$,
but $r_i \neq r_j$, and therefore $\psi$ is violated.
\end{proof}

The following relations play an important role in the relational
description of $\B$ and the clones above $\B$.

\begin{defn}
For $n \geq 2$, let $R_n$, $\underline R_n$, and $R_n^{\neq}$ be the
relations defined by
\begin{align*}
R_n(x_1,y_1,\dots,x_n,y_n) & \equiv \bigvee_{1\leq i \leq n} x_i \neq y_i \\
\underline R_n(x_1,y_1,\dots,x_n,y_n) & \equiv
R_n(x_1,y_1,\dots,x_n,y_n) \vee
x_1=y_1=\dots=x_n=y_n \\
R^{\neq}_n(x_1,y_1,\dots,x_n,y_n) & \equiv
R_n(x_1,y_1,\dots,x_n,y_n) \wedge \bigwedge_{1\leq i \neq j \leq n}
(x_i \neq y_j \wedge x_i \neq x_j \wedge y_i \neq y_j) \, .
\end{align*}
\end{defn}
Note that all three relations can be defined by connected extended
Horn formulas. Observe also that the expressive power of the
relations in each of these sequences increases with increasing $n$:
For example, $\exists u,v.\ (u=v)\wedge
R_{n+1}(x_1,y_1,\ldots,x_n,y_n,u,v)$ is a pp definition of $R_n$
from $R_{n+1}$.

\begin{prop}\label{prop:RelBarChar}
    Let $R$ be a relation with a first-order definition in $({\mathbb N}, =)$. Then the following are equivalent.
    \begin{enumerate}
    \item $R$ is preserved by $\B$, i.e., $\B \subseteq \pol(R)$.
    \item Every expanded Horn formula that defines $R$ is connected extended Horn.
    \item $R$ can be defined by a connected extended Horn formula.
    \item There exists an $n$ such that $R$ has a pp-definition
      in $(\mathbb N,\underline R_n,\neq)$.
    \end{enumerate}
\end{prop}

\begin{proof}
    The implication from $(1)$ to $(2)$ is shown in Lemma~\ref{lem:connextHorn}.
    For the implication from $(2)$ to $(3)$, recall that
    $\B$ contains $\H$ and hence
    the relation $R$ has a Horn definition. By Lemma~\ref{lem:expand},
    every Horn formula has an equivalent expanded Horn formula,
        and by assumption this formula is a connected extended Horn definition of $R$.

    To show that $(3)$ implies $(4)$, let $R$ be a relation with a connected extended Horn definition $\Phi$.
    We first show that
    every extended Horn clause $\psi_1 \rightarrow \psi_2$ from
    $\Phi$ can be pp-defined in $(\mathbb N,\underline R_n)$, for a sufficiently
    large $n$, if $\psi_2$ is non-empty and does not contain \false.

    Let $G$ be the graph of $\psi_1$ where we add isolated
    vertices for each variable that appears in $\psi_2$ but not in $\psi_1$;
    in other words, let $G$ be the graph of $\psi_1 \rightarrow \psi_2$
    where we remove an edge between $x$ and $y$ if $x=y$
    is a literal of $\psi_2$ and not of $\psi_1$.
    Let $C_1,\dots,C_c$ be the components of the graph $G$. We claim that $\psi_1 \rightarrow \psi_2$ can be pp-defined
    by $\underline R_n$ where $n=|C_1|+\dots+|C_c|$.
    A pp definition is obtained from the formula
    \begin{align*}
        \underline R_n(x_1,y_1,\dots,x_n,y_n) \wedge \bigwedge_{(1 \leq i<c)} \bigwedge_{(1 < j  \leq |C_{i+1}|)} y_{s(i)+1} = x_{s(i)+1} = x_{s(i)+j}
    \end{align*}
    by existentially quantifying $x_1,\dots,x_n$, where $s(i) := |C_1|+\dots+|C_{i}|$.
    In this pp definition, the variables $y_{s(i)+1}, \dots,y_{s(i)+|C_{i+1}|}$ are the variables in the component
    $C_{i+1}$ from $\psi_1 \rightarrow \psi_2$.

    To see that this is correct, consider a satisfying
    assignment to the variables of $\psi_1 \rightarrow \psi_2$.
    Set for all $i$ the variables $x_{s(i)+1},\dots,x_{s(i)+|C_{i+1}|}$
    to the same value as $y_{s(i)+1}$.
    Suppose first
    that there exist two variables $y_k,y_l$ from the same component $C_i$
    which are mapped to different values. Then $x_k=x_l$ cannot be
    equal to both $y_k$ and $y_l$, and hence either $x_k \neq y_k$
    or $x_l \neq y_l$. Hence,
    the conjunct $\underline R_n(x_1,y_1,\dots,x_n,y_n)$ is
    satisfied, and the assignment satisfies the given
    pp formula.
    Suppose now that otherwise all variables
    which have the same component are mapped to the same value.
    Since $\psi_1 \rightarrow \psi_2$ is connected Horn, this
    implies by transitivity of equality that all variables that
    appear in $\psi_1 \rightarrow \psi_2$ must have the same value. But this assignment
    clearly satisfies the given pp formula as well, so
    we are done with one direction of the correctness proof.

    For the other direction, consider a satisfying assignment for the
    pp formula. If all variables are set to the same value, then
    $\psi_1 \rightarrow \psi_2$ will be satisfied by this assignment
    as well, because $\psi_2$ does not contain \false.
    Otherwise, $x_k \neq y_k$ for some $1 \leq k \leq n$
    because of the conjunct $\underline R_n(x_1,y_1,\dots,x_n,y_n)$.
    Let $C_{i+1}$ be the component of $y_k$ in $G$.
    Then $y_{s(i)+1}=x_{s(i)+1}=x_k \neq y_k$. But since
    $C_{i+1}$ is a connected component, this implies that
    one of the equalities in $\psi_1$ must be falsified, showing
    that $\psi_1 \rightarrow \psi_2$ is satisfied by the assignment.

    Now, suppose that $\Phi$ contains a clause
    $\psi_1 \rightarrow \psi_2$ where $\psi_2$ contains \false;
    i.e., the clause is logically equivalent to
    a disjunction of literals of the form $x \neq y$.
    We proceed just as in the case above, but
    use the relation $\underline R_{n+1}(x_1,y_1,\dots,x_{n+1},y_{n+1})$
    instead of $\underline R_n(x_1,y_1,\dots,x_{n},y_{n})$ with two new existentially quantified variables $x_{n+1}$ and $y_{n+1}$ at the two additional arguments,
    and we also add the conjunct $x_{n+1} \neq y_{n+1}$ to the formula.
    We leave the verification that the resulting formula is equivalent to $\psi_1 \rightarrow \psi_2$ to the reader.
    As we have seen, any clause from a connected Horn formula can be pp-defined in $(\mathbb N,\underline R_{n+1},\neq)$, and thus
    $(3)$ implies $(4)$.

    Finally, $(4)$ implies $(1)$:
    It is easy to see that for all $n \geq 2$ the relations $\underline R_{n}$ and the relation $\neq$ are preserved by $1$-Bar operations.
\end{proof}

\subsubsection{Above the bar clone}
\label{sssect:abovebar} In this section we show that the interval
between $\B$ and $\pol(\neq)$ contains an infinite strictly
decreasing sequence of clones whose intersection is $\B$.

\begin{defn}\label{def:fk}
    For all $k\geq 2$, fix an operation $f_k \in \O^{(k)}$ satisfying the following:
    \begin{align*}
        f_k(x,1,1,\dots,1,1) & = 1 \\
        f_k(2,x,2,\dots,2,2) & = 2 \\
        & \dots \\
        f_k(k,k,\dots,k,x) & = k
    \end{align*}
    Moreover, for all other arguments, the function arbitrarily takes a value that is distinct from all other function      values.
\end{defn}

Observe that $f_3$ appeared already in
Definition~\ref{def:oddClone}. We use the $f_k$ to define
another sequence of operations.

\begin{defn}\label{def:gk}
For all $k\geq 3$, fix an operation $g_k \in \O^{(k+1)}$ which
satisfies
$$g_k(0,x_1,\dots,x_k) = f_k(x_1,\dots,x_k)$$
and which takes distinct values for all other arguments.
\end{defn}

Observe that each of the operations $g_k$ (and similarly $f_k$)
generates $\B$: Obviously, $g_k$ depends on all arguments, is
non-injective, and preserves $\neq$. By
Proposition~\ref{prop:BarChar}, all relations that are preserved by
$g_k$ are also preserved by $\B$, and hence $g_k$ generates $\B$.

We now show that the operations $g_k$, for increasing $k$, generate
smaller and smaller clones.

\begin{lem}\label{lem:gkpresrk-1}
For $k \geq 3$, the operations $g_k$ and $f_k$ preserve $\underline
R_{k-1}$.
\end{lem}
\begin{proof}
We show the lemma for $g_k$; the proof for $f_k$ is similar and a
bit easier, and left to the reader. Let $t^0,\dots,t^k$ be
$(2k-2)$-tuples that all satisfy $\underline R_{k-1}$. Let
$t:=g_k(t^0,\dots,t^k)$, and suppose that $t_1=t_2, \dots,
t_{2k-3}=t_{2k-2}$. By the definition of $g_k$, this implies that
$t^0=(0,\dots,0)$, and that for all $1 \leq i \leq k-1$ the tuples
$(t^1_{2i-1},\dots,t^{k}_{2i-1})$ and $(t^1_{2i},\dots,t^k_{2i})$
have equal entries except for at most one position, call it $p(i)$.
Since the function $p$ takes at most $k-1$ values, there must be a
$j \in \{1,\dots, k\}$ which is not in the range of $p$. For this
$j$ we thus have $t^j_1=t^j_2, \dots, t^j_{2k-3}=t^j_{2k-2}$.
Because $t^j \in \underline R_{k-1}$, this implies that
$t^j_1=\dots=t^j_{2k-2}$. Hence, the tuple $t=g_k(t^0,\dots,t^k)$ is
constant as well; therefore, $g_k$ preserves $\underline R_{k-1}$.
\end{proof}

\begin{lem}\label{lem:gkviolatesrk}
For $k \geq 3$, the operation $f_k$ preserves neither $\underline
R_k$ nor $R_k$. The operation $g_k$ does not preserve $\underline
R_k$.
\end{lem}
\begin{proof}
If we apply $f_k$ to the $2k$-tuples
\begin{align*}
t^1 & = (0,1,2,2,3,3,\dots,k,k) \\
t^2 & = (1,1,0,2,3,3, \dots, k,k) \\
 & \dots \\
t^k & = (1,1,2,2,3,3,\dots,0,k)
\end{align*}
that all satisfy $\underline R_{k}$ (and all satisfy $R_k$), then we
obtain $(1,1,2,2,3,3,\dots,k,k)$, which does not satisfy $\underline
R_{k}$ (and does not satisfy $R_k$). For $g_k$, we define $t^0 =
(0,\dots,0)$, and apply $g_k$ to $(t^0,t^1,\dots,t^k)$.
\end{proof}

We have thus seen that $\Pol(\underline R_{3},\neq), \Pol(\underline
R_{4},\neq), \dots$ is a strictly decreasing chain of clones (it is
decreasing by the remark after the definition of these relations,
and strictly so by the preceding lemma). By
Proposition~\ref{prop:RelBarChar} its intersection equals $\B$.

\subsubsection{The odd clone}
A clone of special interest (see
Theorems~\ref{thm:monoidalIntervalOfS} and
\ref{thm:monoidalIntervalOfSPlusConst}) is the clone $\S$ generated
by $f_3$, which will be called the \emph{odd clone}
(cf.~Definition~\ref{def:oddClone}). We will now describe the
relations invariant under $\S$ by providing a (finite) generating
system of these relations, as well as a syntactic characterization
of the formulas defining such relations.


\begin{defn}
    Let $\text{ODD}_3$ be the ternary relation $$\{(a,b,c) \in {\mathbb N}^3 \; | \; a=b=c \; \vee |\{a,b,c\}| = 3\} \; .$$
\end{defn}

\begin{prop}\label{prop:SChar}
    Let $R$ be a relation with a first-order definition in $({\mathbb N}, =)$. Then the following are equivalent.
    \begin{enumerate}
    \item $R$ can be defined by a connected Horn formula.
    \item $R$ has a pp definition by $\textit{\underbar R}_2$ and $\neq$.
    \item $R$ has a pp definition by $\text{ODD}_3$ and $\neq$.
    \item $R \in \Inv(\S)$, i.e., $R$ is preserved by $f_3$.
    \item $R \in \Inv(\B)$,
      but there is no definition of $R_3^{\neq}$ in $({\mathbb N},R,\textit{\underbar R}_2,\neq)$.
    \end{enumerate}
\end{prop}

Before proving the proposition, we need two proof-theoretical lemmas
observing restrictions on what can be derived from connected
extended Horn formulas.

\begin{lem}\label{lem:restrSem}
  Suppose $\Phi$ is a connected extended Horn formula that does not imply $a=b$,
  but $\Phi \wedge x=y$ implies $a=b$. Then $\Phi \wedge x=y$
  also implies $x=y=a=b$.
\end{lem}

\begin{proof}
    We
    rewrite $\Phi$ as a conjunction $\Psi$ of (standard) Horn clauses as described above. When for an extended clause $\phi \rightarrow (\psi_1\wedge \dots \wedge \psi_k)$ for some $i \leq k$ we have $\psi_i=\false$, then
    we first remove all positive literals except $\psi_i$ from the extended clause,
    which leads to an equivalent formula.

    It is well-known that the restricted form of resolution for predicate logic that
    requires that at least one of the parent clauses in a resolution
    step has to be positive, also known as \emph{P-resolution},
    is complete  (see e.g.~\cite{SchoeningLogic}, p.~102).
    When we apply resolution on a formula, we assume that
    there are always clauses $u \neq v \vee v \neq w \vee u=w$,
        $u \neq v \vee v = u$, and $u=u$ for universally quantified
        variables $u,v,w$ (that axiomatize that the symbol $=$ for equality is an equivalence relation), and we also assume that
    all other variables in the formula are in fact constant symbols (i.e., they are existential variables before skolemization).

    Since $P$-resolution is complete, we have in particular that
    $\Psi \wedge x=y$ implies $a=b$ if and only if $a=b$ can be derived by P-resolution from $\Psi \wedge x=y$.
    We show that $\Psi \wedge x=y$ implies
    $x=y=a=b$ for all formulas
    $a=b$ that are not implied by
    $\Psi$ but implied by $\Psi \wedge x=y$.
    Suppose for contradiction
    that there is a literal $a=b$ that can be derived by P-resolution from  $\Psi \wedge x=y$ but cannot be derived from $\Psi$, and where
    $\Psi \wedge x=y$ does not imply $x=y=a=b$.
    Select a literal $a=b$ with these properties
    whose $P$-derivation from $\Psi \wedge x=y$ is shortest possible.

    The literal $a=b$ has either been derived from a clause $a_1 \neq b_1 \vee \dots \vee a_k \neq b_k \vee a=b$ in $\Psi$ by resolving with the literals $a_1 = b_1, \dots, a_k = b_k$, or by
    resolving with of the three clauses having universally quantified
    variables (after instantiating the universal variables).
    In the second case, suppose that the universal clause was
    $u \neq v \vee v \neq w \vee u=v$, and $a=b$
    has been derived by resolving with $a=c$
    and $c=b$, for some $c$ (for the other two universal clauses, the
    argument is analogous).
    Then at least one of the literals $a=c$ and $c=b$
    cannot be derived from $\Psi$ alone,
     because otherwise $a=b$ can be derived from $\Psi$ alone,
  in contradiction to our assumptions. Suppose the
    literal that cannot be derived from $\Psi$ alone is $a=c$
    (the other case is analogous).
    By the choice of $a=b$, the formula
    $\Psi \wedge x=y$ implies $x=y=a=c$.
    At the same time, we have that $\Psi \wedge x=y$ implies
    $a=b$, and hence $\Psi$ implies $x=y=a=b$, a contradiction.

    In the first case, at least one of the literals in the clause
    $a_1 \neq b_1 \vee \dots \vee a_k \neq b_k \vee a=b$,
    say $a_i=b_i$, cannot be derived from $\Psi$ alone,
  because otherwise $a=b$ can be derived from $\Psi$ alone,
  in contradiction to our assumptions.
    By the choice of $a=b$, $\Psi \wedge x=y$ implies $x=y=a_i=b_i$.
    Let $\phi$ be the extended Horn clause in $\Phi$ that
    corresponds to the Horn clause $\psi$ in $\Psi$
    (i.e., $\psi$
    is among the clauses that have replaced
    the clause $\phi$ when creating $\Psi$ from $\Phi$).

    Since $\phi$ is a connected extended Horn clause, there is a sequence of variables $x_1,x_2,\dots,x_s$ of $\phi$ with $x_1=a$ and $x_s=a_i$ such that $x_i=x_{i+1}$ for $1 \leq i < s$
    is either a sub-formula of the right-hand side or of the left-hand
    side of the implication $\phi$. By construction of $\Psi$,
    for each literal $x_i=x_{i+1}$
    on the right-hand side there is the clause
    $a_1 \neq b_1 \vee \dots \vee a_k \neq b_k \vee x_i = x_{i+1}$
    in $\Psi$. But then it is clear that $\Psi \wedge a_1 = b_1 \wedge \dots \wedge a_k = b_k$ implies that $x_1 = \dots = x_s$,
    and hence also that $a_i=b_i=a=b$. Therefore, $\Psi \wedge x=y$
    implies $x=y=a=b$, a contradiction.
\end{proof}

The following is a straightforward consequence.
\begin{lem}\label{lem:restrSem2}
 Let $\Phi$ be a connected extended Horn formula
 and let $\psi$ be a conjunction of equalities $x_1=y_1 \wedge \dots \wedge x_n = y_n$.
 If $\Phi \wedge \psi$ implies $x_0=y_0$, but $\Phi \wedge \psi'$
 does not imply $x_0=y_0$ for any subconjunction $\psi'$ of $\psi$,
  then $\Phi \wedge \psi$
  also implies $x_0=y_0=x_1=y_1=\dots=x_n=y_n$.
\end{lem}
\begin{proof}
    We apply Lemma~\ref{lem:restrSem}
    to $\Phi' := \Phi \wedge x_2=y_2 \wedge \dots \wedge x_n=y_n$, $x_1=y_1$, and $x_0 = y_0$
    to derive that $\Phi'  \wedge x_1=y_1$
    implies $x_0=y_0=x_1=y_1$. We can argue analogously
    for $x_i=y_i$ instead of $x_1=y_1$ for all $1 \leq i \leq n$,
    and the statement follows.
\end{proof}

\begin{proof}[Proof of Proposition~\ref{prop:SChar}]
    To show that (1) implies (2), suppose that $\Phi$
    is a connected Horn definition of $R$.
    We show that each
    clause $\psi$ of $\Phi$
    has a pp definition in $(\mathbb N,\textit{\underbar R}_2,\neq)$.

    First, consider the case that $\psi$
    contains a positive literal $x_0=y_0$.
    The graph of $\psi$ without the edge between $x_0$ and $y_0$ has at
    most two connected components. In the case that the graph is still
    connected even after removing the edge between $x_0$ and $y_0$, the
    clause $\psi$ is a tautology and can be ignored.
    Otherwise, let $x_0,\dots,x_l$ be the variables of the
    connected component of $x_0$, and let $y_0,\dots,y_k$ be the variables
    of the connected component of $y$.

    Consider the formula $\psi'$ defined by
    \begin{align*}
        \exists u_1,\dots,u_l,v_1,\dots,v_k. \;
        \textit{\underbar R}_2(x_0,x_1,u_1,u_1) \wedge \textit{\underbar R}_2(u_1,x_2,u_2,u_2) \wedge \dots \wedge \textit{\underbar R}_2(u_{l-1},x_l,u_l,u_l) \\
        \wedge \, \textit{\underbar R}_2(y_0,y_1,v_1,v_1) \wedge \textit{\underbar R}_2(v_1,y_2,v_2,v_2) \wedge \dots \wedge \textit{\underbar R}_2(v_{k-1},y_k,v_k,v_k)
        \wedge u_l = v_k \; .
    \end{align*}
    We show that a tuple satisfies $\psi$ if and only if it satisfies $\psi'$.
    If the tuple satisfies no negative literal of $\psi$, then
    the only situation where $\psi'$ is true is
    $x_0=u_1=\dots=u_l=v_k=\dots=v_1=y_0$; hence in this case, $\psi'$ holds if and only if
    the tuple satisfies $x_0=y_0$, which in turn is the case if and only if $\psi$ holds. If on the other hand the tuple
    satisfies some negative literal of $\psi$, then either there is $i<l$ such that $x_i \neq x_{i+1}$ or $i<k$ such that $y_i\neq y_{i+1}$; assume wlog the
    first is the case. Then we set $(u_1,\dots,u_i)$ to $(x_1, \dots, x_i)$
    and $(v_1,\dots,v_k)$ to $(y_1,\dots,y_k)$, $u_l$ to $v_k$, and
    set $u_{i+1},\dots,u_{l-1}$ to arbitrary pairwise disjoint values that
    are also distinct from all the values taken by any other variable
    of the clause. It is easy to check that this assignment of variables
    is a satisfying assignment for the existentially quantified variables in
    the formula $\psi'$. Hence, $\psi'$ is a pp definition of
    $\psi$.

    If $\psi$ does not contain a positive literal,
    then the graph of $\psi$
    consists of at most two connected components, because $\phi$ is connected Horn.
    Assume first it does have two components. Let $x_0,\dots,x_l$ be
    the variables in the first and $y_0,\dots,y_k$ be the variables
    in the second connected component.
    As above, it is easy to verify
    that the following formula is equivalent to $\psi$.
   \small
    \begin{align*}
      \exists u_1,\dots,u_l,v_1,\dots,v_k,z. &
      \textit{\underbar R}_2(x_0,x_1,u_1,u_1) \wedge \textit{\underbar R}_2(u_1,x_2,u_2,u_2) \wedge
      \dots \wedge \textit{\underbar R}_2(u_{l-1},x_l,u_l,u_l) \\
    &
      \wedge \, \textit{\underbar R}_2(y_0,y_1,v_1,v_1) \wedge \textit{\underbar R}_2(v_1,y_2,v_2,v_2) \wedge
      \dots \wedge \textit{\underbar R}_2(v_{k-1},y_k,v_k,v_k) \\
      & \wedge \textit{\underbar R}_2(x_l,u_l,z,z) \wedge \textit{\underbar R}_2(y_k,v_k,z,z)
      \wedge u_l \neq v_k \; .
    \end{align*}
    \normalsize

    If the graph of $\psi$ has only one component, then the
    following formula is equivalent to $\psi$:
    \small
    \begin{align*}
      \exists u_1,\dots,u_l&
      \textit{\underbar R}_2(x_0,x_1,u_1,u_1) \wedge \textit{\underbar R}_2(u_1,x_2,u_2,u_2) \wedge
      \dots \wedge \textit{\underbar R}_2(u_{l-1},x_l,u_l,u_l)\wedge u_1\neq u_l \; .
    \end{align*}
    \normalsize

    For the implication from $(2)$ to $(3)$,
    note that  $\textit{\underbar R}_2(x_1,x_2,x_3,x_4)$ can be pp defined by
    \begin{align*}
     \exists y_1,y_2,y_3,y_4,y_5. & \; \text{ODD}_3(x_1,x_2,y_1) \wedge \text{ODD}_3(x_1,x_2,y_2) \\
    & \wedge \text{ODD}_3(y_1,y_2,y_3) \wedge \text{ODD}_3(y_3,y_4,y_5) \\
    & \wedge \text{ODD}_3(x_3,x_4,y_4) \wedge \text{ODD}_3(x_3,x_4,y_5)
    \end{align*}
    We check that this formula is indeed equivalent to $\textit{\underbar
    R}_2(x_1,x_2,x_3,x_4)$:
    Let  $t$ be a $4$-tupel that satisfies
    $\textit{\underbar R}_2(x_1,x_2,x_3,x_4)$; we assign values to
    the $y_i$ such that the formula above is satisfied.
    If $t_1=t_2$ then set $y_1, y_2, y_3$ to $t_1$.
    If $t_3=t_4$ then set $y_3,y_4,y_5$ to $t_3$.
    Note that if both $t_1=t_2$ and $t_3=t_4$ then
    our choice for $y_3$ is well-defined, since
    $\textit{\underbar R}_2(x_1,x_2,x_3,x_4)$ then implies that $t_1=t_2=t_3=t_4$.
    If $t_1 \neq t_2$ then set $y_1,y_2$ to distinct values
    that are also distinct from all entries of $t$.
    Similarly, if $t_3 \neq t_4$, then set $y_4,y_5$ to distinct values that
    are also distinct from all entries of $t$.
    If both $t_1 \neq t_2$ and $t_3 \neq t_4$, then set $y_3$
    to any value that is distinct from the values for $y_1,y_2,y_4,y_5$.
    It is straightforward to verify that these assignments satisfy
    all conjunct in the given formula. Conversely, if
    a $4$-tuple $t$ satisfies the formula, and $t_1=t_2$ and $t_3=t_4$,
    then $t_1=t_2=t_3=t_4$, and hence $t$ satisfies $\textit{\underbar R}_2(x_1,x_2,x_3,x_4)$.

    For the implication from $(3)$ to $(4)$
    it suffices to verify that $f_3$ preserves $\text{ODD}_3$.
    Every 3-tuple $t$ from $\text{ODD}_3$ either satisfies $t_1=t_2=t_3$ or $t_1
    \neq t_2 \neq t_3 \neq t_1$. To violate $\text{ODD}_3$, we would have to find three
    such tuples $a,b,c$ such that $d:=f_3(a,b,c)$ has
    two equal entries that are distinct from a third entry, say wlog.\
    $d_1=d_2 \neq d_3$. The only two ways to have
    $d_1=d_2$ are when $d_1=d_2$ is from
    $\{0,1,2\}$, or when this is not the case and the tuples $(a_1,b_1,c_1)$ and $(a_2,b_2,c_2)$ are
    equal. Consider the first case, and assume wlog.\ $d_1=d_2=0$. Then we
    know that $b_1=c_1=0$ and $b_2=c_2=0$. But then, since $b$ and $c$
    are in $\text{ODD}_3$, we have $b_3=c_3=0$ as well, and hence
    $d_1=d_2=d_3=0$, a contradiction. The second case is similar and left to the reader.
     Hence, $f_3$ preserves $\text{ODD}_3$.

    For the implication from $(4)$ to $(5)$, containment
    in $\Inv(\B)$ follows from the fact that $f_3$ generates all operations in $\B$ (see also the remark after Definition~\ref{def:gk}).
    By Lemma~\ref{lem:gkpresrk-1},
    $f_3$ preserves $\textit{\underbar R}_2$ and $f_3$ also clearly preserves
    $\neq$, so it suffices to show that $f_3$
    violates $R^{\neq}_3$.
    If we apply $f_3$ to the 6-tuples $t_1 = (0,1,2,2,3,3)$,
    $t_2 = (0,0,1,2,3,3)$, $t_3 = (0,0,2,2,1,3)$, which all satisfy
    $R^{\neq}_3$, we obtain $(0,0,2,2,3,3)$, which does not
    satisfy $R^{\neq}_3$.

  Finally, we show the contraposition of the implication from $(5)$ to $(1)$: if $R$ cannot be defined by a connected Horn formula,
  but is preserved by $\B$,
  then there is a pp-definition of $R^{\neq}_n$
  in $({\mathbb N},R,\textit{\underbar R}_2,\neq)$ for some $n \geq 3$.
  Because $R^{\neq}_3$ clearly has a pp-definition by
  $R^{\neq}_n$ for all $n \geq 3$ this will show the statement.

  Let $\Phi(v_1,\dots,v_s)$ be an expanded Horn formula that
  defines $R$. By Proposition~\ref{prop:RelBarChar}, we
  know that $\Phi$ is a connected extended Horn formula. We first
  claim that
  if for all extended Horn clauses $\psi_1 \rightarrow \psi_2$
  in $\Phi$ the graph of $\psi_1$ has
  at most two connected components, then
  $\Phi$ has a connected Horn definition, contradicting our
  assumption on $\Phi$. To see this, we show that every such clause $\psi_1 \rightarrow
  \psi_2$ has a connected Horn definition.
  If $\psi_2$ contains \false\ or is empty then there is nothing to
  show; assume thus that this is not the case.
    If the graph of $\psi_1$ has exactly two components
  $C_1$ and $C_2$, we know that $\psi_2$ contains a conjunct $x=y$ where $x \in C_1$ and $y \in C_2$, because $\Phi$ is connected
  extended Horn.
  Then $\psi_1 \rightarrow \psi_2$
  is equivalent to the conjunction of the connected Horn clauses $\psi_1 \rightarrow x=y$
  and $(\psi_1 \wedge x=y) \rightarrow x'=y'$ for every
  other conjunct $x'=y'$ of $\psi_2$.   The case that the graph of
  $\psi_1$ is connected is easier and left to the reader.

 Therefore, $\Phi$ must contain a clause
 $\psi_1 \rightarrow \psi_2$ such that $\psi_1$ consists of $n \geq 3$ connected components $C_1,\dots,C_n$.
 We show that $R^{\neq}_n(x_1,y_1,\dots,x_n,y_n)$
has the following pp-definition by $R$, $\textit{\underbar R}_2$ and
$\neq$.

\begin{align*}
\exists v_1,\dots,v_s. & \; R(v_1,\dots,v_s) \wedge \bigwedge_{1 \leq i \leq n} \bigwedge_{v_k \in C_i} \textit{\underbar R}_2(x_i,y_i,v_k,v_k) \\
& \wedge \bigwedge_{1\leq i \neq j \leq n} (x_i \neq y_j \wedge x_i
\neq x_j \wedge y_i \neq y_j) \nonumber && (*)
\end{align*}

 To prove the equivalence, assume that $t \notin
R^{\neq}_n$. We want to show that $t$ does not satisfy $(*)$. If $t$
violates the second line of the formula $(*)$, we are done.
Otherwise, it must hold that $x_1=y_1,\dots,x_n=y_n$. Suppose for
contradiction that there are elements $v_1,\dots,v_s$ as specified
in $(*)$. Then for all $1 \leq i \leq n$ the sub-formula
$\bigwedge_{v_k \in C_i} \textit{\underbar R}_2(x_i,y_i,v_k,v_k) $
implies that all variables $v_k$ from $C_i$ must have the same value
as $x_i=y_i$. Then the clause $\psi_1 \rightarrow \psi_2$ in $\Phi$
implies that $x_1=y_1= \dots = x_n = y_n$ (because $\psi_1
\rightarrow \psi_2$ is connected), which contradicts the second line
of formula $(*)$.

Now assume that $t \in R^{\neq}_n$. We want to show that $t$
satisfies $(*)$. If $t$ satisfies $x_i \neq y_i$ for all $1 \leq i
\leq n$, then any tuple from $R$ gives a assignment to the
existentially quantified variables in $(*)$ that shows that $t$
satisfies $(*)$. For the sake of notation, let us thus assume that
there is a $k$ with $1 \leq k < n$ such that $t$ satisfies $x_i=y_i$
for $1 \leq i \leq k$, and $x_i \neq y_i$ for all $k<i\leq n$.
We claim that there exists a tuple $r \in R$ such that
 for all $i,j \leq k$ and all $v_l \in C_i$, $v_m \in C_j$ it holds that $r_l=r_m$
 if and only if $i=j$.
Otherwise, the formula
$$\Phi(v_1,\dots,v_s) \wedge \bigwedge_{i
\leq k} \bigwedge_{v_l,v_m \in C_i} v_l = v_m \wedge \bigwedge_{v_l
\in C_i, v_m \in C_j, 1\leq i\neq j \leq k} v_l \neq v_m$$ is
unsatisfiable. Since P-resolution (also see the proof of
Lemma~\ref{lem:restrSem}) is complete, there is a P-resolution proof
of unsatisfiability. This proof uses at most once a clause of the
form $\{v_p \neq v_q\}$ from this formula, since it resolves only
with positive literals and since therefore every resolution step
using $\{v_p \neq v_q\}$ must be the last one, deriving
\textit{false}. Therefore, $\Phi(v_1,\dots,v_s) \wedge \bigwedge_{i
\leq k} \bigwedge_{v_l,v_m \in C_i} v_l = v_m$ implies $v_p=v_q$.
Let $\phi$ be the smallest sub-formula of $\bigwedge_{i \leq k}
\bigwedge_{v_l,v_m \in C_i} v_l = v_m$ such that $\Phi \wedge \phi$
implies $v_p=v_q$. If $\phi$ is empty, then the clause $v_p=v_q$ is
contained in $\Phi$ because $\Phi$ is expanded; but then, adding
$v_p\neq v_q$ to $\psi$ in the clause $\psi_1 \rightarrow \psi_2$
results in an equivalent formula, a contradiction to $\Phi$ being
expanded. If $\phi$ is not empty, we can apply
Lemma~\ref{lem:restrSem2} and obtain that $\Phi \wedge \phi$ implies
that all variables $u_1,\dots,u_l$ that appear in $\phi$ have the
same value as the value of $v_p=v_q$.
Because $\Phi$ is expanded, it has to contain the connected extended
Horn clause $\phi \rightarrow u_1 = \dots = u_l = v_p=v_q$; again,
this shows that adding $v_p\neq v_q$ to $\psi_1$ in the clause
$\psi_1 \rightarrow \psi_2$ results in an equivalent formula, a
contradiction.

Hence, we have shown the existence of a tuple $r \in R$ with the
properties as claimed above. The tuple $r$ gives witnesses for the
existentially quantified variables in $(*)$ that show that $t$
satisfies $(*)$.
\end{proof}
    \subsection{Richard}
        \label{subsect:richard}

In this section we show that $\R$ and $\S$ are incomparable clones,
and that every clone in $I_\I$ is either contained in $\S$ or
contains $\R$. We start with the fundamental observation that the
set of operations $\R$ (see
Definition~\ref{def:injectiveInOneDirectionAndRichard}) is indeed a
clone.

\begin{lem}\label{lem:RisAClone}
    $\R$ is a clone.
\end{lem}

\begin{proof}
    Clearly, $\R$ contains the projections. Let $f\in\R\un$ and
    $g_1,\ldots,g_n\in\R\um$. There exist $1\leq i\leq n$ and
    $1\leq j\leq m$ such that $f$ is injective in the $i$-th
    direction and $g_i$ is injective in the $j$-th direction. It
    is readily verified that $f(g_1,\ldots,g_n)$ is injective in
    the $j$-th direction.
\end{proof}

\begin{defn}
Let $\phi$ be a quantifier-free first-order formula where all atomic
subformulas are of the form $x=y$. Then $\phi$ is called
\emph{negative} if $\phi$ is in conjunctive normal form and
    each clause in $\phi$ either consists of a single positive literal,
    or does not contain positive literals.
\end{defn}

The following theorem has been stated in equivalent form
in~\cite[Lemma 8.6]{qecsps}. For the convenience of the reader, and
for completeness, we also present a simplified proof here.

\begin{thm}\label{thm:negative}
     Let $R$ be a relation having a reduced Horn definition that is not negative. Then $\text{ODD}_3$ has a pp definition from $R$ and $\neq$.
\end{thm}

\begin{proof}
Let $\Phi(x_1,\dots,x_n)$ be a reduced Horn definition of $R$ which
is not negative; assume moreover that among such definitions of $R$,
$\Phi$ is one with the minimal possible number of occurrences of
variables. Since $R$ is not negative, there must be a clause $\psi$
in $\Phi$ of the form $x_{i_1} = x_{j_1} \vee x_{i_2} \neq x_{j_2}
\vee \dots \vee x_{i_k} \neq x_{j_k}$ for $k \geq 2$. Since $\Phi$
is reduced, it contains a tuple $a \in \mathbb N^n$ with $a_{i_l} =
a_{j_l}$ for all $1 \leq l \leq k$: Otherwise, $\Phi$ would imply
$x_{i_1} \neq x_{j_1} \vee x_{i_2} \neq x_{j_2} \vee \dots \vee
x_{i_k} \neq x_{j_k}$, and we would obtain an equivalent formula by
removing the literal $x_{i_1} = x_{j_1}$ from $\psi$. Moreover, $R$
contains a tuple  $b \in \mathbb N^n$ with $b_{i_1} \neq b_{j_1}$,
$b_{i_2} \neq b_{j_2}$, and $b_{i_l} = b_{j_l}$ for all $3 \leq l
\leq k$: Otherwise, we could remove the literal $x_{i_2} \neq
x_{j_2}$ from $\psi$.

 Now consider the formula
\begin{align*}
R(x_1,\dots,x_n) \wedge \bigwedge^k_{l \geq 3} x_{i_l} = x_{j_l} &&
(*)
\end{align*}
Note that because $R$ is reduced, we have $|\{i_1,j_1,i_2,j_2\}|
\geq 3$.

We first look at the case where $|\{i_1,j_1,i_2,j_2\}| = 3$. By
appropriately renaming the variables of $\Phi$, we may assume that
$i_1=i_2=1$, and that $j_1=2$ and $j_2=3$. Let
$R'(x_i,x_{j_1},x_{j_2})$ be the ternary relation defined by
existentially quantifying all variables in (*) except for
$x_i$,$x_{j_1}$,$x_{j_2}$. Because $a \in R$ satisfies
$a_{i_l}=a_{j_l}$ for all $3\leq l \leq k$, it satisfies (*). Since
$a_{j_1}=a_{i_1}=a_{i_2}=a_{j_2}$, the relation $R'$ contains a
triple with three equal entries. Moreover, because $b \in R$
satisfies (*) too, the relation $R'$ contains a tuple
$(b_1,b_2,b_3)$ with $b_1 \neq b_2$ and $b_1 \neq b_3$.

Assume first that $R'$ contains a triple with three distinct
entries. We know that $R'$ does not contain any tuple of the form
$(x,x,y)$ for $x \neq y$, because every $n$-tuple $a$ with
$a_{i_1}=a_{j_1}$ and $a_{i_2} \neq a_{j_2}$ violates $\psi \wedge
 \bigwedge_{l=3}^k x_{i_l} = x_{j_l}$.
Hence, the relation defined by $R'(x,y,z) \wedge R'(y,z,x) \wedge
R'(z,x,y)$ only contains the tuples with three distinct entries and
the tuples with three equal entries, and we have indeed found a
pp-definition of the relation $\text{ODD}_3$ in $(\mathbb N,R)$.

Striving for a contradiction, we suppose now that
$R'(x_{i},x_{j_1},x_{j_2})$ does not contain a triple with three
distinct entries. Then it is impossible that there is a tuple $r \in
R'$ where $r_{2} \neq r_3$: This is because having a Horn
definition, $R$ is preserved by a binary injection, by
Proposition~\ref{prop:HornChar}, and so is $R'$. But the tuple
obtained by applying a binary injective operation to $r$ and
$(b_1,b_2,b_3)$ has three distinct entries, which contradicts our
assumption for this case. Thus the clause $x_{i_3} \neq x_{j_3} \vee
\dots \vee x_{i_k} \neq x_{j_k} \vee x_{j_1}=x_{j_2}$ is entailed by
$\Phi$. Hence, we could replace the clause $\psi$ in $\Phi$ by this
shorter clause and obtain a formula that is equivalent to $\Phi$,
because conversely, $x_{j_1}=x_{j_2}$ implies $x_{i_1}=x_{j_1} \vee
x_{i_1} \neq x_{j_2}$ (recall that $i_1=i_2$). This contradicts the
choice of $\Phi$.

Now take the case where $|\{i_1,j_1,i_2,j_2\}| = 4$; say wlog
$i_1=1, j_1=2,i_2=3,j_2=4$. Let
$R'(x_{i_1},x_{j_1},x_{i_2},x_{j_2})$ be the relation defined by
$(*)$ by existentially quantifying all variables except for
$x_{i_1}, x_{j_1}, x_{i_2}, x_{j_2}$. As before we observe that $a'
:= (a_{i_1},a_{j_1},a_{i_2},a_{j_2}) \in R'$ and $b' :=
(b_{i_1},b_{j_1},b_{i_2},b_{j_2}) \in R'$. If $a_{i_1}=a_{j_1} \neq
a_{i_2}=a_{j_2}$, then the tuple obtained by applying a binary
injective operation to $a'$ and $b'$ has four distinct entries; this
tuple is an element of $R$, by Proposition~\ref{prop:HornChar}. It
is easy to verify that then the formula
\begin{align*} R'(x_1,x_2,x_3,x_4) \wedge R'(x_3,x_4,x_1,x_2) \wedge x_1 \neq x_3 \wedge x_1 \neq x_4 \wedge x_2 \neq x_3 \wedge x_2 \neq x_4 && (**)
\end{align*}
is a pp definition of the relation $T(x_1,x_2,x_3,x_4)$ given by
$$T := \{(x_1,x_2,x_3,x_4) \; | \; x_1=x_2 \neq x_3=x_4  \vee \, |\{x_1,x_2,x_3,x_4\}|=4 \} \, .$$
The relation $\text{ODD}_3(x_1,x_2,x_3)$ has the following pp-definition.
$$\exists y \exists z. \; T(x_1, x_2, y, z) \wedge T(x_2, x_3, y, z) \wedge T(x_1,x_3,y,z) \; .$$

Now suppose that $a_{j_1} = a_{i_2}$ and hence $a'$ is the tuple
with four equal entries. Assume for a moment that there exist $p \in
\{1,2\}$ and $q \in \{3,4\}$ such that $R'(x_1,\dots,x_4)$ implies
$x_p = x_q$. Then the relation obtained from $R'(x_1,\dots,x_4)$ by
existentially quantifying $x_q$ is such that we can proceed as in
the case where $|\{i_1,j_1,i_2,j_2\}|=3$. Thus, we may henceforth
assume that $x_p = x_q$ is not entailed for any $p \in \{1,2\}$ and
$q \in \{3,4\}$. For every such pair $(p,q)$, fix a tuple in $R'$
witnessing that $x_p \neq x_q$ can be satisfied. Now applying a
$5$-ary injection to these four witnessing tuples and $b'$, we
obtain a tuple all of whose entries are distinct. This tuple is in
$R'$ as $R'$ is preserved by $\H$, by
Proposition~\ref{prop:HornChar}.

Consider the formula $R''(x_1,x_2,x_3,x_4)$ defined by the following
formula.

\begin{align*}
\bigwedge_{\pi \in S_{\{1,2,3,4\}}}
R'(x_{\pi(1)},x_{\pi(2)},x_{\pi(3)},x_{\pi(4)}),
\end{align*}
where $S_{\{1,2,3,4\}}$ denotes the symmetric group on
$\{1,2,3,4\}$. We claim that $R''$ defines $\text{ODD}_3$. Note that tuples of
the form $(x,x,x,y)$ and $(x,x,y,z)$, for distinct $x,y,z$, are not
contained in $R'$ since every $n$-tuple $d$ with $d_{i_1}=d_{j_1}$
and $d_{i_2} \neq d_{j_2}$ violates $\psi \wedge
 \bigwedge_{l=3}^k x_{i_l} = x_{j_l}$.
Therefore, the relation $R''$ does not contain tuples with two
values where one value appears three times, or tuples with three
values where one value appears twice. But certainly $R''$ does
contain all tuples with four equal values, and all tuples with four
pairwise distinct values. If some tuple where two values occur twice
is excluded from $R'$, then $R''$ does not contain any tuple with
two values. Hence in this case, $\exists u. R''(x_1,x_2,x_3,u)$ is a
pp definition of $\text{ODD}_3(x_1,x_2,x_3)$, and we are done. Otherwise, $R''$
contains all tuples where two values occur twice. Similarly as
before it is easy to verify that then the expression~$(**)$ above
where $R'$ has been replaced by $R''$ is a pp definition of
$T(x_1,x_2,x_3,x_4)$. We have already seen that $\text{ODD}_3$ has a pp
definition from $T$, so it has a pp definition in $(\mathbb N, R,
\neq)$ as well.

\end{proof}

The following proposition describes the clone $\R$: Item (1) shows
that $\R$ is finitely generated (namely, by any binary operation
which violates $\text{ODD}_3$ but is injective in one direction); and (2) and
(3) provide a syntactical description of the formulas defining
relations in $\Inv(\R)$.
\begin{prop}\label{prop:disjofineq}
    Let $R$ be a relation with a first-order definition in $({\mathbb N}, =)$. Then the following are equivalent.
    \begin{enumerate}
    \item $R$ is preserved by a binary
    operation from $\pol(\neq)$ that violates $\text{ODD}_3$.
    \item Every reduced definition of $R$ is negative.
    \item $R$ has a negative definition.
    \item $R$ is preserved by $\R$.
    \end{enumerate}
\end{prop}

\begin{proof}
    To see that $(1)$ implies $(2)$, fix any reduced definition of $R$. Since $R$ is preserved by a binary
    operation from $\pol(\neq)$, Proposition~\ref{prop:HornChar}
    implies that this definition is Horn. As there is a binary
    operation preserving $R$ and $\neq$ which violates $\text{ODD}_3$, there is no pp
    definition of $\text{ODD}_3$ from $R$ and $\neq$. Hence, our reduced
    definition must be negative, by Theorem~\ref{thm:negative}.

    The implication from $(2)$ to $(3)$ is trivial.

    To show that $(3)$ implies $(4)$, fix a negative definition $\Phi$
    of $R$. Let $f\in\R\un$ and $a_1,\ldots,a_n\in R$ be arbitrary.
    Assume without loss of generality that
    $f$ is injective in the $1$-st direction. Then
    $b=f(a_1,\ldots,a_n)$ satisfies all inequalities that were
    satisfied by $a_1$, which is enough to satisfy $\Phi$ as
    $\Phi$ is negative. Hence, $f$ preserves $R$.

    The implication from $(4)$ to $(1)$ is witnessed by any binary
    operation which is not in $\B$ but injective in one direction, as is easily verified.
\end{proof}

The following proposition establishes statement (4) of
Theorem~\ref{thm:monoidalIntervalOfS}.

\begin{prop}\label{prop:rsAntichain}
    The clones $\R$ and $\S$ form a maximal antichain in the interval
    $[\B,\pol(\neq)]$. In fact, every clone in this interval is either contained in $\S$ or contains
    $\R$. In particular, $\R$ is a cover of $\R\cap\S$.
\end{prop}

\begin{proof}
First, we have to show that $\R$ and $\S$ are indeed incomparable.
By Proposition~\ref{prop:disjofineq} it is clear that $\R$ is not a
subset of $\S$, because it contains an operation that violates $\text{ODD}_3$.
To show that $\S$ is not a subset of $\R$, it suffices to show that
the operation $f_3$ from Definition~\ref{def:oddClone} is not
injective in any direction. Indeed, this holds as $f_3(x,1,1)=1$ and
$f_3(2,y,2)=2$ and $f_3(3,3,z)=3$ for all $x,y,z\in \mathbb N$.
Therefore, $f\nin\R$ and hence $\R$ and $\S$ form an antichain.

Now, let $\C$ be any clone from the interval $[\B,\pol(\neq)]$ which
is not contained in $\S$. Then it contains an operation in
$\pol(\neq)$ that does not preserve $\text{ODD}_3$. Since $\text{ODD}_3$ has two orbits,
$\C$ contains a binary operation with these properties, by
Lemma~\ref{lem:op-arity}. Thus, Proposition~\ref{prop:disjofineq}
implies that $\C$ is above $\R$.

Hence, the antichain that consists of $\R$ and $\S$ is maximal, and
every clone of the interval not contained in $\S$ does contain $\R$.
\end{proof}
    \subsection{Richard's many friends}
        \label{subsect:richardsmanyfriends}

We now prove that there exists an antitone and injective mapping
from the power set of $\omega$ into the interval $[\R,\pol(\neq)]$;
in particular, this interval has cardinality $2^{\aleph_0}$. This
proves the last statement (5) of
Theorem~\ref{thm:monoidalIntervalOfS}.

\begin{defn}
    For all $n\geq 3$, write
    $$
        \delta_n:=x_1\neq y_1\vee\ldots\vee x_n\neq y_n.
    $$
    For all $A\subseteq\{1,\ldots,n\}$ with $1<|A|<n$, writing
    $A=\{j_1,\ldots,j_r\}$ with $j_1<j_2<\ldots<j_r$, we set
    $$
        \kappa_A:= y_{j_1}\neq x_{j_2}\vee y_{j_2}\neq x_{j_3}\vee\ldots\vee y_{j_r}\neq x_{j_1}.
    $$

    Set
    $$
        \gamma_n:=\delta_n\wedge\bigwedge_{A\subseteq\{1,\ldots,n\},1<|A|<n}\kappa_A.
    $$
\end{defn}

Observe that the $\gamma_n$ are all negative formulas; hence, by
Proposition~\ref{prop:disjofineq}, polymorphism clones of sets of
relations defined by such $\gamma_n$ will always contain $\R$. In
the rest of this section, we will prove that no fixed $\gamma_n$ can
be defined from the others by primitive positive definitions.
Therefore, distinct sets of such formulas define distinct clones
above $\R$, which is what we want.

\begin{defn}
    We now enumerate tuples $c\in \nat^{2n}$ as
    $c=(c^{1,x},c^{1,y},\ldots,c^{n,x},c^{n,y})$. If
    $\phi(x^1,y^1,\ldots,x^n,y^n)$ is a formula, then we say that $c$ satisfies
    $\phi$ iff $\phi(c^{1,x},c^{1,y},\ldots,c^{n,x},c^{n,y})$
    holds; this is to say that we insert $c^{i,u}$ for $u^i$, where $u\in\{x,y\}$ and $1\leq i\leq n$. Moreover, if $c_1,\ldots,c_m\in \nat^{2n}$, we denote the ``$(i,x)$-th column tuple'' $(c^{i,x}_1,\ldots,c^{i,x}_m)$ by
    $c^{i,x}$.
\end{defn}

\begin{defn}
    For all $n\geq 3$, we set $C_n$ to be the $2n$-ary relation defined by $\gamma_n$.
\end{defn}

The operations of the following definition will separate the
polymorphism clones of the $C_n$.

\begin{defn}
    For all $n\geq 3$, the \emph{Hubie-violator} $H_n$ is an
    operation defined as follows: Enumerate
    $C_n\cap\{1,\ldots,n+1\}^{2n}$, that is, those elements of $C_n$
    which have only entries in $\{1,\ldots,n+1\}$, by
    $c_1,\ldots,c_m$. Note that $m>0$ since $C_n$ always contains the tuple where $x_i=0$ and $y_i=1$, for all $i \leq n$.
    Now we define $H_n\in\Om$ by $H_n(c^{j,x})=H_n(c^{j,y})=j$,
    for all $1\leq j\leq n$; in other words, we set $H_n(c_1,\ldots,c_m)=(1,1,\ldots,n,n)$.
    For every other input tuple in $\nat^m$, $H_n$ takes a unique value distinct from all other values in its range.\\
\end{defn}

\begin{lem}
    Let $n\geq 3$. Then the Hubie-violator $H_n$ violates $C_n$.
\end{lem}
\begin{proof}
    By definition, $H_n(c_1,\ldots,c_m)=(1,1,\ldots,n,n)$, a tuple which does not satisfy $\delta_n$ and which
    is therefore not an element of
    $C_n$.\\
\end{proof}

\begin{lem}
    Let $3\leq k<n$. Then the Hubie-violator $H_n$ preserves $C_k$.
\end{lem}
\begin{proof}
    Write $m$ for the arity of $H_n$. Let $t_1,\ldots,t_m\in C_k$ and set $s=H_n(t_1,\ldots,t_m)$. We must show $s\in C_k$.

    We first check that $s$ satisfies the clauses $\kappa_A$. Suppose it does
    not, and pick any $A\subseteq \{1,\ldots,k\}$ with $1<|A|<k$ witnessing this. We may assume
    wlog that $A=\{1,\ldots,p\}$, where $1<p<k$, and write
    $s=(a_p,a_1,a_1,a_2,a_2,a_3,\ldots
    ,a_{p-1},a_{p-1},a_p,?,\ldots)$. Set $Q:=\{a_j:1\leq j\leq
    p\}$, and $Q_1:=Q\cap\{1,\ldots,n\}$, and $Q_2=Q\setminus Q_1$.

    If $Q_1$ were empty,
    then since for every value in $Q_2$ there is only a unique tuple that $H_n$ sends to that value, we
    have that the column
    $t^{i,y}$ equals the column $t^{i+1,x}$, for all $1\leq i\leq p$ (where
    we set $p+1:=1$). Thus in that case,
    no row $t_j$ would satisfy $\kappa_A$, a contradiction.

    Assuming
    that $Q_2$ were empty, we now show that there exists $1\leq j\leq m$ such that the row $t_j$ does not
    satisfy $\kappa_{\{1,2\}}$. Since $H_n(t^{1,x})=a_p$, we have $t^{1,x}=c^{a_p,u}$, for some $u\in\{x,y\}$.
    Similarly, $t^{1,y}=c^{a_1,w}$, $t^{2,x}=c^{a_1,q}$,
    and
    $t^{2,y}=c^{a_2,v}$, where
    $w,q,v\in\{x,y\}$. Consider first the case where
    $a_p\neq a_1$ and $a_2\neq a_1$ ($a_p$ may well equal $a_2$).
    Pick $d\in C_n\cap\{1,\ldots,n+1\}^{2n}$ such that
    $d^{a_1,x}=d^{a_1,y}=a_1$, $d^{a_p,u}=a_2$, and $d^{a_2,v}=a_2$. This is
    without doubt possible by defining all other entries to be very
    unequal. Because in the definition of the Hubie-violator all tuples of $C_n\cap\{1,\ldots,n+1\}^{2n}$ appear,
    there is $1\leq j\leq m$ with $d=c_j$. We have:
    $t_j^{1,x}=c_j^{a_p,u}=a_2$,  $t_j^{1,y}=c_j^{a_1,w}=a_1$,
     $t_j^{2,x}=c_j^{a_1,q}=a_1$, and $t_j^{2,y}=c_j^{a_2,v}=a_2$. Hence, $t_j=(a_2,a_1,a_1,a_2,?,\ldots)$
    does not satisfy $\kappa_{\{1,2\}}$, a contradiction. The case where
    $a_p=a_1$ or $a_2=a_1$ is even easier: Say wlog $a_p=a_1$, and pick $c_j\in C_n$ such that
    $c_j^{a_1,x}=c_j^{a_1,y}=a_1$, and $c_j^{a_2,v}=a_1$. We have:
    $t_j^{1,x}=c_j^{a_1,u}=a_1$,  $t_j^{1,y}=c_j^{a_1,w}=a_1$,
     $t_j^{2,x}=c_j^{a_1,q}=a_1$, and $t_j^{2,y}=c_j^{a_2,v}=a_1$. Hence, in this case
    $t_j=(a_1,a_1,a_1,a_1,?,\ldots)$ does not
    satisfy $\kappa_{\{1,2\}}$ either, a contradiction.

    Finally,
    consider the case where neither $Q_1$ nor $Q_2$ are empty.
    Assume wlog that $s^{1,x}=a_p\in Q_1$ and $s^{1,y}=a_1 \in Q_2$. We have $t^{1,x}=c^{a_p,v}$, for some $v\in\{x,y\}$. Let $r\leq p$
    be minimal with the property that $s^{r,y}\in Q_1$; write $s^{r,y}=a_i$.
    Then $t^{r,y}=c^{a_i,u}$ for some $u\in\{x,y\}$. Pick $c_j\in C_n$
    such that $c_j^{a_i,u}=c_j^{a_p,v}=a_1$. We then have that $t_j$
    does not satisfy $\kappa_{\{1,\ldots,r\}}$: Indeed $t_j^{1,x}=a_1=t_j^{r,y}$, and $t^{h,y}=t^{h+1,x}$ for all $1\leq h< r$ as columns since
    $a_h\in Q_2$. A contradiction.

    It remains to show that $s$ satisfies $\delta_k$. Suppose it does not; we claim that under this assumption, there is $1\leq j\leq m$
such that $t_j$ does not satisfy $\delta_k$ either. Set $Q$ to
    consist of all values that appear in $s$, $Q_1:=Q\cap\{1,\ldots,n\}$, and $Q_2:=Q\setminus
    Q_1$.  Write $p$ for the cardinality of $Q_1$; wlog we may assume $Q_1=\{1,\ldots,p\}$. Observe that since $s$ has length $2k$ and since it does not satisfy $\delta_k$, we have $p\leq k$. Consider the $2n$-tuple
    $b=(1,1,2,2,\ldots,p,p,p+1,n+1,p+2,n+1,\ldots,n,n+1)$. Using $p\leq k<n$, one
    readily checks that $b\in C_n$, hence $b=c_j$ for some $1\leq j\leq m$. Thus, since $c_j$ appears in the $j$-th row of the definition of the Hubie-violator $H_n$, we have that for all
    $d\in \nat^m$ and all $1\leq i\leq p$, if $H_n(d)=i$, then $d_j=i$.
    In particular, this holds for the vectors $t^{r,x}$ and
    $t^{r,y}$ for which $H_n(t^{r,x})=H_n(t^{r,y})\in Q_1$, and we have $t^{r,x}_j=t^{r,y}_j$ for those columns. But for the other columns $t^{r,x},t^{r,y}$ whose (equal) image under $H_n$ is in $Q_2$, we have $t^{r,x}_j=t^{r,y}_j$ anyway since values outside $\{1,\ldots,n\}$ are taken by at most one tuple. Summarizing these two observations, we have that
    $t_j^{r,x}=t_j^{r,y}$ for all $1\leq r\leq k$, so $t_j$
    does not satisfy $\delta_k$. A contradiction.
\end{proof}

We thus know that for all $n\geq 3$,  $C_n$ has no pp-definition
from $\{C_3,\ldots, C_{n-1}\}$. For otherwise, the polymorphisms of
$\{C_3,\ldots,C_{n-1}\}$ would be contained in those of $C_n$, which
is not the case by the preceding two lemmas. It remains to show that
$C_n$ cannot be defined from the $C_k$ with $k>n$ either.

\begin{lem}
    Let $3\leq k<n$. Then the Hubie-violator $H_k$ preserves $C_n$.
\end{lem}
\begin{proof}
    Denote the arity of $H_k$ by $m$. As before, let $t_1,\ldots,t_m\in C_n$ and set $s:=H_k(t_1,\ldots,t_m)$. To see that $s$
    satisfies the clauses $\kappa_A$, one proceeds as in the preceding proof, with $k$ and $n$ exchanged.
    What is more difficult here is to check that $s$ satisfies $\delta_n$. Suppose it does not, and define $Q,Q_1,Q_2$, as well as $p=|Q_1|$,
    as before (only with $k$ and
    $n$ switched). Now we distinguish two cases, $p<k$ and $p=k$.

    In the first
    case, we argue like at the end of the proof of the preceding lemma: Consider the tuple
    $b=(1,1,2,2,\ldots,p,p,p+1,k+1,p+2,k+1,\ldots,k,k+1)$. One
    readily checks that $b\in C_k$, so $b=c_j$ for some $1\leq j\leq m$. Thus, since $c_j$ appears in the $j$-th line of the definition of $H_k$, we have that for all
    $d\in \nat^m$ and all $1\leq i\leq p$, if $H_k(d)=i$, then $d_j=i$.
    In particular, this holds for those vectors $t^{r,x}$ and
    $t^{r,y}$ which are sent by $H_k$ to a value in $Q_1$. Since $t^{r,x}=t^{r,y}$ for the other columns (i.e., those with value in $Q_2$), this implies that
    $t_j^{r,x}=t_j^{r,y}$ for all $1\leq r\leq n$, so $t_j$
    does not satisfy $\delta_n$, a contradiction.

    To finish the proof, consider the case where $p=k$.  Observe first that for every $1\leq i\leq k$, there exists exactly one $1\leq r\leq n$ such that $s^{r,x}=s^{r,y}=i$. There is at least one such $r$ as $p=k$; suppose there were distinct $r,r'$ with this property. Then $t^{r,x},t^{r,y},t^{r',x},t^{r',y}$ are all elements of $\{c^{i,x},c^{i,y}\}$. Pick any tuple $c_j\in C_k$ with $c^{i,x}=c^{i,y}$. We have that $t^{r,x}_j,t^{r,y}_j,t^{r',x}_j,t^{r',y}_j$ are all equal, hence $t_j$ does not satisfy $\kappa_{\{r,r'\}}$, a contradiction.
We may therefore wlog assume that $s$ looks like this:
$s=(1,1,\ldots,k,k,?,\ldots)$, where all entries starting from the
question mark are elements of $Q_2$.  If we had $t^{i,x}=t^{i,y}$
for some $1\leq i \leq k$, then we could derive a contradiction just
like in the case $p<k$, by taking $i$ out of $Q_1$. So assume this
is not the case. Then
    the set of vectors
    $\{t^{1,x},t^{1,y},\ldots,t^{k,x},t^{k,y}\}$ equals
    $\{c^{1,x},c^{1,y},\ldots,c^{k,x},c^{k,y}\}$.
We prove that there exists $1\leq j\leq m$ such that
    $t_j$ does not satisfy $\kappa_{\{1,\ldots,k\}}$. Let $\sigma$ be the permutation on the indices of the tuple $(c^{1,x},c^{1,y},\ldots,c^{k,x},c^{k,y})$ which sends this tuple
    to $(t^{1,x},t^{1,y},\ldots,t^{k,x},t^{k,y})$. Observe that $\sigma$ only switches $c^{j,x}$ and $c^{j,y}$, if
    necessary. Now consider the tuple
    $b:=(k,1,1,2,2,\ldots,k-1,k-1,k)\in C_k$, and apply
    $\sigma^{-1}$ to the indices of this tuple, obtaining a new tuple $\sigma^{-1}(b)$. One
    readily checks $\sigma^{-1}(b)\in C_k$, so $\sigma^{-1}(b)=c_j$ for some $1\leq j\leq
    m$. Thus the tuple $(t_j^{1,x},t_j^{1,y},\ldots,t_j^{k,x},t_j^{k,y})$, which is obtained by applying $\sigma$ to $c_j=\sigma^{-1}(b)$, equals
    $b$. But $b$ does not satisfy
    $\kappa_{\{1,\ldots,k\}}$, so $t_j\nin C_n$, a contradiction.
\end{proof}

We end this section by the following proposition, which completes
the proof of Theorem~\ref{thm:monoidalIntervalOfS} and of the
corollary after.

\begin{prop}
    The mapping $\sigma$ from the power set of $\omega$ into $\Cl_{loc}(\mathbb{N})$ which sends every $A\subseteq \omega$ to
    the local clone $\pol(\{C_{n+3}:n\in A\})$ is injective and antitone (with respect to inclusion, on both sides). In particular, the number
    of local clones containing $\sn$ equals the continuum.
\end{prop}
\begin{proof}
     Since $\pol$ is antitone, it follows that the same holds for
     $\sigma$. Now let $A, B\subseteq \omega$ be unequal, say wlog that there is $n\in
     A\setminus B$. Then the Hubie-violator $H_{n+3}$ violates
     $C_{n+3}$ but preserves all $C_{k+3}$, where $k\in B$. Thus
     $H_{n+3}\in\sigma(B)\setminus\sigma(A)$, and $\sigma$ is indeed
     injective.
\end{proof}

\section{Clones with essential infinite range operations plus constants}\label{sect:infiniteRangeOpsPlusConstants}

It remains to investigate one last monoidal interval, namely the one
corresponding to the monoid $\I^+$ consisting of all unary
operations which are either injective or constant. We will thus be
concerned with the proof of
Theorem~\ref{thm:monoidalIntervalOfSPlusConst} in this section. The
crucial ingredients to the proof will be
Proposition~\ref{prop:rsAntichain} and the following lemma; the
definition of $\S$ was given in Definition~\ref{def:oddClone}.

\begin{lem}\label{lem:S+locallyclosed}
    $\cl{\S^+}=\S^+$.
\end{lem}
\begin{proof}
    We have to show that $\S^+$ is a local clone. Let $f\in\cl{\S^+}$ be an arbitrary non-constant operation. We claim that $f\in\S$. To see this, let $A$ be any finite    subset of $\nat^n$, where $n$ is the arity of $f$. Extend $A$ to a finite set $B\subseteq \nat^n$ such that $f$ is non-constant on $B$. Since $f\in\cl{\S^+}$, and since $\S$ is by definition generated by $f_3$, there exists a term $t$ using $f_3$ and the constants which agrees with $f$ on $B$. Since $t$ is non-constant, it can as well be written without constants, using $f_3$ only: This is because inserting a constant as an argument of $f_3$ gives us an injection or essentially a bar operation, both of which are generated by $f_3$ without using constants. Thus, $f$ can be interpolated on $B$, and hence also on $A$, by a term in $\S$. This proves $f\in\S$.
\end{proof}

Recall that $\R$ consists of all operations which are injective in
one direction
(Definition~\ref{def:injectiveInOneDirectionAndRichard}). Whereas
$\S^+$ is a local clone, we cannot add constants to $\R$ without
generating all operations:

\begin{prop}
    $\cl{\R^+}=\O$.
\end{prop}
\begin{proof}
    Let $f(x_1,\ldots,x_n)\in\O$ be arbitrary, and let $A\subseteq \nat^n$ be finite. Define $g(x_1,\ldots,x_{n+1})$ as follows: If $x_{n+1}=0$ and $(x_1,\ldots,x_n)\in A$, then $g$ returns $f(x_1,\ldots,x_n)$. For every other input tuple, $g$ returns a unique value distinct from all other values in its range. Clearly, $g$ is injective in the $(n+1)$-st direction, so $g\in\R$. Moreover, $g(x_1,\ldots,x_n,0)$ agrees with $f$ on $A$, implying $f\in\cl{\R^+}$ as  $A$ was arbitrary.
\end{proof}

Note that the preceding proposition proves item (7) of
Theorem~\ref{thm:monoidalIntervalOfSPlusConst}. We have found the
largest element of our monoidal interval $I_{\I^+}$:

\begin{lem}\label{lem:S+isTheLargest}
    $\pol(\I^+)=\S^+$.
\end{lem}
\begin{proof}
    Clearly, the unary fragment of $\S^+$ equals $\I^+$, hence $\S^+$ is an element of the corresponding monoidal interval $I_{\I^+}$. Suppose there was a local clone $\C$ in this interval which properly contains $\S^+$. Take any $f\in\C\setminus\S^+$. Clearly, $f$ must be essential. Since the only non-injections in the unary fragment of $\C$ are constant, Lemma~\ref{lem:s5:HaddadRosenberg:forInfinite} implies that $f$ must have infinite range. By Proposition~\ref{prop:s5:essentialWithInfiniteImagePreservesNeq}, $f$ preserves $\neq$. Thus, $f\in\pol(\neq)\setminus\S$, so $f$ generates $\R$ by Proposition~\ref{prop:rsAntichain}. Hence $\C=\O$ by the preceding proposition.
\end{proof}

The following proposition proves statement (1) of
Theorem~\ref{thm:monoidalIntervalOfSPlusConst}.

\begin{prop}
    Let $\C\in I_{\I^+}$. Then $\C^-=\C\cap\S$. In particular, $\C^-$ is a local clone in
        $I_\I$.
\end{prop}
\begin{proof}
    By Lemma~\ref{lem:S+isTheLargest}, $\C$ is a subclone of $\S^+$; hence intersecting $\C$ with $\S$ just means taking away the constants, which proves our assertion. Since the intersection of local clones is a local clone, so is $\C^-$. Clearly, the unary fragment of $\C^-$ equals $\I$, hence $\C^-$ is an element of $I_\I$.
\end{proof}

Observe that the second statement of the theorem is trivial. We turn
to the proof of item~(3).

\begin{prop}
The mapping from $I_{\I^+}$ into the subinterval $[\cl{\I},\S]$ of
$I_\I$  which sends every clone $\C$ to $\C^-$ is a complete lattice
embedding which preserves the smallest and the largest element.
\end{prop}
\begin{proof}
    It is clear that this mapping preserves the order and by the second statement of Theorem~\ref{thm:monoidalIntervalOfSPlusConst}, it sends the smallest (largest) element of $I_{\I^+}$ to the smallest (largest) element of  $[\cl{\I},\S]$. It is also obvious that the mapping is injective. Let $\C,\D\in I_{\I^+}$. Clearly, $(\C\cap\D)^-=\C^-\cap\D^-$, hence the mapping preserves finite meets; larger meets work the same way. Since $(\C\vee\D)^-$ contains both $\C^-$ and $\D^-$, it contains also their join, so $(\C\vee\D)^-\supseteq \C^-\vee\D^-$. For the other inclusion, let $f$  be any non-constant operation in the (not necessarily local) clone generated by $\C\cup\D$. Then $f$ can be written as a term over $\C\cup\D$. We can assume that in this term, thinking of it as a tree, the constants appear only as leaves. We may also assume that except for leaves, no subtree of this tree is constant. That just means that the nodes above the leaves are operation symbols in $\C^-\cup\D^-$, some of whose variables are set to constant values. But these operations in $\C^-\cup\D^-$ with some variables set to constants are again elements of $\C^-\cup\D^-$, as they are non-constant by assumption and since $\C$ and $\D$ are clones containing all constants. Summarizing, we can write $f$ as a term over $\C^-\cup\D^-$ without the use of constants. Hence, $f\in \C^-\vee\D^-$. This proves $(\C\vee\D)^-\subseteq \C^-\vee\D^-$, and hence our mapping preserves finite joins. Infinite joins work the same way.
\end{proof}

Observe that item (4) of
Theorem~\ref{thm:monoidalIntervalOfSPlusConst} is clear from the
definitions. We prove (5).

\begin{prop}
    For all $\C\in I_\I$ which do not contain $\R$, $\cl{\C^+}$
        is a local clone in $I_{\I^+}$.
\end{prop}
\begin{proof}
    By Proposition~\ref{prop:rsAntichain}, $\C$ is contained in $\S$. Hence, $\cl{\C^+}$ is contained in $\S^+$ by Lemma~\ref{lem:S+locallyclosed}, which equals $\pol(\I^+)$ by Lemma~\ref{lem:S+isTheLargest}. Therefore $\cl{\C^+}$ is indeed a local clone in $I_{\I^+}$.
\end{proof}

Item (6) is trivial from the definitions. We finish the proof of the
theorem by the following proposition, which restates items (8) and
(9):

\begin{prop}
    $\H^+$ is the unique cover of $\cl{\I^+}$ and  $\B^+$  is the unique cover of $\H^+$  in $I_{\I^+}$.
\end{prop}
\begin{proof}
    Let $\C$ be in $I_{\I^+}$, and assume it has an essential operation $f$. Then $f\in\S$ by Lemma~\ref{lem:S+isTheLargest}, thus $f$ generates $\H$ by Proposition~\ref{prop:HornChar}. It is easy to see that $\H^+$ is a local clone, which proves the first assertion. Similarly, if $\C$ is in $I_{\I^+}$ and properly contains $\H^+$, then any operation witnessing this generates $\B$, by Proposition~\ref{prop:BarChar}. Again, it is straightforward to check that $\B^+$ is a local clone.
\end{proof}

We want to close this section with the remark that the mapping which
sends every clone $\C$ in the interval $[\cl{\I},\S]$ of $I_\I$ to
$\cl{\C^+}$ is not injective. In particular, it collapses all clones
between $\R \cap \S$ and $\S$.

\begin{prop}\label{prop:+notinjective}
$\cl{(\R \cap \S)^+} = \S^+$
\end{prop}
\begin{proof}
The operation $g_3$ (Definition~\ref{def:gk}) is injective in the
first argument, and hence contained in $\R$. By Lemma 59, $g_3$
preserves $\underline R_2$, and, being an element of $\R$, it also
preserves $\neq$. Proposition~\ref{prop:SChar} then implies that
$g_3$ is contained in $\S$: Otherwise, it would violate a relation
preserved by $\S$; that relation would have a pp definition from
$\underline R_2$ and $\neq$, a contradiction. Therefore, the
operation $g_3$ is in $\R \cap \S$.

The operation given by $g_3(0,x_1,x_2,x_3)$ equals by definition the
operation $f_3$ (Definition~\ref{def:oddClone} or \ref{def:fk}),
which by definition generates $\S$. Hence, $\cl{(\R \cap \S)^+}$
contains $\S^+$. The converse containment is trivial.
\end{proof}

Observe that this proposition implies that the complete lattice embedding that sends every $\C\in I_{\I^+}$ to $\C^-$ is not surjective onto the interval  $[\cl{\I},\S]$: If there existed a clone $\C\in I_{\I^+}$ such that $\C^-=\R\cap \S$, then $\C=(\R\cap\S)^+$ since $\C=(\C^-)^+$. Thus $(\R\cap\S)^+$ would be a clone and hence equal to $\S^+$, implying $\R\cap\S=\S$, a contradiction. In fact the same argument shows that a clone $\D\in [\cl{\I},\S]$ is in the range of this embedding iff $\D^+$ is a clone.

\section{Open Problems}
We have to leave the complete description of the monoidal intervals corresponding to $\I$ and $\I^+$ open; we do not know if a reasonable description of these intervals is possible at all. Two particular problems seem most important to us at this point:
\begin{itemize}
    \item Can one effectively decide whether
for a given sequence $R_0,R_1,\dots,R_n$ of relations that are first-order definable in $(\mathbb N,=)$
there is a pp-definition of $R_0$ in the structure $(\mathbb N, R_1,\dots,R_n)$?
    \item What is the cardinality of the monoidal interval corresponding to $\I^+$?
\end{itemize}

\vspace{.4cm}

\paragraph{Acknowledgements}
We thank the referee for his helpful suggestions.

\end{document}